\documentclass[reqno, 11pt]{article}
\usepackage{amssymb, amsmath, amsthm, amsfonts, amscd, epsfig}
\usepackage{color}

\usepackage[english]{babel}
\usepackage{bbm}
\usepackage{graphicx, epsfig}
\usepackage{geometry,graphicx,pgfplots,caption}

\usepackage[export]{adjustbox}
\usepackage[titletoc,toc,title]{appendix}
\usepackage[normalem]{ulem} 

\newtheorem{theorem}{Theorem}[section]
\newtheorem{cor}[theorem]{Corollary}
\newtheorem{lemma}[theorem]{Lemma}                                                                                                                                                                                                                                                                             
\newtheorem{prop}[theorem]{Proposition}
\newtheorem{definition}{Definition}

\newtheorem{remark}{Remark}
{\bf}{\it}

\def\ep{\epsilon}

\newcommand{\Om}{\Omega}

\newcommand{\RR}{\mathbb{R}}
\newcommand{\CC}{\mathbb{C}}
\newcommand{\NN}{\mathbb{N}}
\newcommand{\ZZ}{\mathbb{Z}}
\newcommand{\Scal}{\mathcal{S}}

\newcommand{\Kcal}{\mathcal{K}}

\newcommand{\teta}{{\widetilde{\eta}}}
\newcommand{\tomega}{{\widetilde{\omega}}}
\newcommand{\p}{\partial}

\newcommand{\ds}{\displaystyle}
\newcommand{\eqnref}[1]{(\ref {#1})}

\newcommand{\beq}{\begin{equation}}
\newcommand{\eeq}{\end{equation}}
\newcommand{\RN}[1]{%
\textup{\uppercase\expandafter{\romannumeral#1}}%
}

\newcommand{\eps}{\varepsilon}
\newcommand{\hKcal}{{\hat{\mathcal{K}}}}

\numberwithin{equation}{section}
\numberwithin{figure}{section}

\begin{document}

\newcommand{\TheTitle}{The Eigenvalue Problem for the Laplacian via Conformal Mapping and the Gohberg--Sigal Theory}

\newcommand{\TheAuthors}

\author{
Marius Beceanu\thanks{Guest Researcher, Stephen Redenti Research Group, Lehman College, CUNY (mbeceanu@gmail.com).}
\and Jiho Hong\thanks{Postdoctoral Fellow, Department of Mathematics, CUHK, Shatin, N.T., Hong Kong (jihohong@cuhk.edu.hk).}
\and Hyun-Kyoung Kwon\thanks{Department of Mathematics and Statistics, University at Albany, SUNY, 1400 Washington Ave., Albany, NY, 12222, USA (hkwon6@albany.edu)}
\and Mikyoung Lim\thanks{Department of Mathematical Sciences, Korea Advanced Institute of Science and Technology, 291 Daehak-ro, Yuseong-gu, Daejeon 34141, Republic of Korea (mklim@kaist.ac.kr).}}
\title{{\TheTitle}\thanks{{This research was supported by the Basic Science Research Program through National Research Foundation of Korea (NRF) grants NRF-2019R1A6A1A10073887, NRF-2021R1A2C1011804 and RS-2024-00359109.}}}

\maketitle

\begin{abstract}
	We consider the Dirichlet and Neumann eigenvalues of the Laplacian for a planar, simply connected domain. The eigenvalues admit a characterization in terms of a layer potential of the Helmholtz equation. Using the exterior conformal mapping associated with the given domain, we reformulate the layer potential as an infinite-dimensional matrix. Based on this matrix representation, we develop a finite section approach for approximating the Laplacian eigenvalues and provide a convergence analysis by applying the Gohberg--Sigal theory for operator-valued functions. Moreover, we derive an asymptotic formula for the Laplacian eigenvalues on deformed domains that results from the changes in the conformal mapping coefficients.
\end{abstract}

\noindent {\footnotesize {\bf Mathematics subject classification
(MSC2020): 30C35, 35J05, 35P15, 41A25, 47A56} }



%

\section{Introduction}

Let $\Omega$ be a bounded simply connected domain in $\mathbb{R}^2$ with $C^{1,\alpha}$ boundary for some $\alpha>0$.
We first consider the Neumann eigenvalue problem for the Laplacian on $\Omega$ given by
\begin{equation}\label{problem:neumann:zero}\begin{cases}
		\ds -\Delta u=\lambda u\quad\mbox{in }{\Omega},\\
		\ds\frac{\partial u}{\partial \nu}=0\quad\mbox{on }\partial\Omega,
\end{cases}\end{equation}
where $\lambda$ is a constant and $\nu$ is the outward unit normal vector to $\partial \Omega$.
If \eqnref{problem:neumann:zero} admits a nontrivial solution $u$, then we call $\lambda$ a \emph{Neumann eigenvalue} of $-\Delta$.
The Neumann eigenvalues $\lambda_j^{\rm N}$ form a sequence of non-negative real numbers that diverges to infinity,
$$\label{Neumann:increasing}
0=\lambda_1^{\rm N}<\lambda_2^{\rm N}\le\lambda_3^{\rm N}\le\cdots\le\lambda_n^{\rm N}\to\infty.
$$
We similarly define the \emph{Dirichlet eigenvalues} $\lambda_j^{\rm D}$ by replacing the boundary condition in \eqnref{problem:neumann:zero} with the zero boundary condition $u=0$ on $\partial \Omega$. We again enumerate them as 
$$\label{Dirichlet:increasing}
0<\lambda_1^{\rm D}\le\lambda_2^{\rm D}\le\lambda_3^{\rm D}\le\cdots\le\lambda_n^{\rm D}\to\infty.
$$
It is known that for all $j \in \mathbb{N}$, $\lambda_{j+1}^{\rm N}\le\lambda_j^{\rm D}$ \cite{Friedlander:1991:IDN}. Numerous estimates for the eigenvalues of the Laplacian in two dimensions have been investigated as shown in the work \cite{Faber:1923:BDU, Girouard:2009:MSN, Goldshtein:2016:EFV, Henrot:2006:EPE, Krahn:1925:URM, Polya:1961:EVM, Szego:1954:ICE, Weinberger:1956:IIN}. Since the Rayleigh-Ritz method readily gives an upper bound, the emphasis has been on finding a lower bound. In this paper, we are concerned with approximating the eigenvalues of the Laplacian with the Dirichlet or the Neumann boundary condition.

For $\varphi\in L^2(\p\Om)$ and $z \in \mathbb{R}^2 \backslash \partial \Omega$, the \emph{single-} and \emph{double-layer potentials} for the Helmholtz equation are defined by 
\beq
\label{SD}
\begin{aligned}
	\mathcal{S}_\Omega^\omega[\varphi](z)&=\int_{\p\Omega}{\Gamma^\omega(z-\xi)}\varphi(\xi)\,d\sigma(\xi), \\ 
	\mathcal{D}_\Omega^\omega[\varphi](z)&=\int_{\p\Omega}\frac{\p\Gamma^\omega(z-\xi)}{\p\nu_{\xi}}\,\varphi(\xi)\,d\sigma(\xi),
\end{aligned}
\eeq
where $\Gamma^\omega$ denotes the fundamental solution to the Helmholtz equation $(\Delta+\omega^2)u=0$ in two dimensions. Both of these integral operators solve the Helmholtz equation in $\mathbb{R}^2 \backslash \partial \Omega$. For $\omega>0$, $\Gamma^\omega$ admits an absolutely convergent series expansion
\beq\label{eq:fundamental:sol}
\Gamma^\omega(z)
=\frac{1}{2\pi}\log|z|+\tau(\omega) 
+\sum_{j=1}^\infty b_j\big(\log \left(\omega|z|\big) +c_j
\right) \big(\omega|z|\big)^{2j}, \quad z \in \mathbb{R}^2 \backslash \{0\},\eeq
where $\tau(\omega)=\frac{1}{2\pi}\left(\log \frac{\omega}{2} +\gamma_0\right)-\frac{{ \mathrm{i}}}{4}$, $\gamma_0$ stands for the Euler constant, and 
\beq
\label{bjcj}
b_j=\frac{(-1)^j}{2\pi}\frac{1}{{2^{2j}}(j!)^2},\quad
c_j= \gamma_0-\log 2-\frac{\pi {\mathrm{i}}}{2}-\sum_{n=1}^j\frac{1}{n}.
\eeq

We define related integral operators for $\varphi \in L^2(\partial \Omega)$ and $z \in \partial \Omega$ with weakly singular kernels as
\beq
\begin{aligned} \label{Kdef}
	\ds\mathcal{K}_\Om^\omega[\varphi](z)&= \mbox{p.v.} \int_{\p\Omega}\frac{\p\Gamma^\omega(z-\xi)}{\p\nu_{\xi}}\,\varphi(\xi)\,d\sigma(\xi),\\
	\mathcal{K}_\Omega^{\omega,*}[\varphi](z)&=\mbox{p.v.}\int_{\p\Omega}\frac{\p\Gamma^\omega(z-\xi)}{\p\nu_{z}}\,\varphi(\xi)\,d\sigma(\xi),
\end{aligned}
\eeq
where p.v. stands for the Cauchy principal value. As in $\Gamma^{\omega}(z)$, the operator-valued function $\Kcal_\Om^\omega$ can be analytically extended to $\CC\setminus\sqrt{-1}\RR^-,$ where $\RR^-=\{x \in \mathbb{R}: x \leq 0\}$. In particular, it admits a local power series expansion in $\omega$. The relationship between these operators and the previous ones are given by
\beq 
\begin{aligned}\label{relationS} \notag
	\ds \mathcal{D}_{\Omega}^{\omega}[\varphi] \text{ }\vline_{\pm} (z)&=\left( \mp \frac{1}{2}I+\mathcal{K}_{\Omega}^{\omega} \right)[\varphi](z), \\
	\frac{\partial \mathcal{S}^{\omega}_{\Omega}}{\partial \nu}\text{ }\vline_{\pm}(z)&= \left( \pm \frac{1}{2}
	I+\mathcal{K}^{\omega*}_{\Omega} \right) [\varphi](z).
\end{aligned}
\eeq

For a  $C^{1,\alpha}$ domain $\Omega$, the operator $\mathcal{K}_\Omega^0$ is compact on $L^2(\p\Om)$ \cite{Specht:1969:CNP,Taylor:2000:TFP}. 
Since for a fixed $\omega$, $\mathcal{K}_\Omega^\omega-\mathcal{K}_\Omega^0$ is a Hilbert--Schmidt integral operator, it follows that
$\Kcal_\Om^\omega$ is compact on $L^2(\p\Om)$ as well.
Moreover, the operator $\pm\frac{1}{2}I+\mathcal{K}_\Om^\omega$ is Fredholm with index one and is invertible on $L^2(\p\Om)$ for each $\omega\in \CC\setminus\sqrt{-1}\RR^-$ except at most on a discrete set. From the result by M. Mitrea \cite{Mitrea:1996:BHH} (see Theorem \ref{Lemma:eigen:DN} in Subsection \ref{subsec:multiple}), which connects the eigenvalue characterization problem to the problem about the invertibility of these operators, $\omega_0^2$ is a Neumann eigenvalue if and only if $\omega_0$ is 
a \emph{characteristic value} of $-\frac{1}{2}I+\mathcal{K}^{\omega}_{\Omega}$ (the definition of a characteristic value will be presented in the next section). The same holds for a Dirichlet eigenvalue $\omega_0^2$ with $\frac{1}{2}I+\mathcal{K}_{\Omega}^{\omega}$.

Conformal mapping techniques have been powerful tools for dealing with various planar inclusion problems in two dimensions. One particular example stands out; for any simply connected domain, there uniquely exists an exterior conformal mapping that is represented by certain canonical Laurent series (see \eqnref{eq:Psi:laurent} below). A basis $\{\widetilde{\eta}_m\}_{-\infty}^{\infty}$ for $L^2(\partial \Omega)$ previously studied by Y. Jung and the fourth author using the conformal mapping (\cite{Jung:2021:SEL}) contributes to a matrix formulation for $\mathcal{K}^{\omega}_{\Omega}$. Next, the conformal mapping identifies the shape of the domain so that the shape perturbation of the domain can be interpreted as small modifications in the corresponding conformal mapping. 
We note that $\Kcal_\Om^\omega$ is an analytic operator-valued function with respect to $\omega$.
Operator analogues of the Argument Principle and Rouch\'e's Theorem of complex function theory due to  I. C. Gohberg and E. I. Sigal (\cite{Gohberg:1971:OEL}) are then invoked to asymptotically estimate the eigenvalues using the characteristic values of approximations of $\pm \frac{1}{2}I+\mathcal{K}_{\Omega}^{\omega}$. 
More precisely, we consider the following two eigenvalue problems for the Laplacian:
\begin{itemize}
	\item[\rm(i)] approximation of the eigenvalues by the finite section method for  $\mathcal{K}^{\omega}_{\Omega}$ via the geometric basis $\{\widetilde{\eta}_m\}_{-\infty}^{\infty}$.
	
	\item[\rm(ii)] estimation of the effect of small perturbations on the Laplacian eigenvalues.
\end{itemize}

Our main theorems for problem (i) are as follows.
Results concerning the Dirichlet eigenvalues can be stated analogously.
Here, $T_n$ stands for an $n$-th truncation that will be defined precisely later on, and $P_n$ denotes the projection from $L^2(\partial \Omega)$ onto $\text{span }\{\widetilde{\eta}_m:-n \leq m \leq n\}$:

\begin{theorem}\label{theorem:main:first}
	Let $\p\Om$ be of class $C^{1,\alpha}$ for some $0<\alpha<1$. For a Neumann eigenvalue $\omega_0^2 \neq 0$ of $-\Delta$ on $\Omega$ of multiplicity $\mu$, the following hold: 
	\begin{itemize}
		\item[\rm (a)]
		There exists $N\in\NN$ such that for $n\geq N$, $-\frac{1}{2}I+T_n\mathcal{K}_{\Omega}^{\omega}$ on $L^2(\p\Om)$ has, counting multiplicity, exactly $\mu$ characteristic values in some closed ball $\overline{B_{\delta}(\omega_0)}$ for constant $\delta>0$. 
		\item[\rm (b)] There exists a constant $C>0$ such that
		$$
		\left|\omega_{n,j}-\omega_0\right|\le\frac{C}{n!},\quad\mbox{for all }n\geq N,\ j=1,\dots,\mu,$$
	\end{itemize}
	where the $\omega_{n,j}$ ($j=1,\dots,\mu$) denote the characteristic values of $-\frac{1}{2}I+T_n\mathcal{K}_{\Omega}^\omega$. 
\end{theorem}

\begin{theorem}\label{theorem:main:second}
	Let $\p\Om$ be of class $C^{1,\alpha}$ for some $0<\alpha<1$. For a Neumann eigenvalue $\omega_0^2 \neq 0$ of $-\Delta$ on $\Omega$ of multiplicity $\mu$, the following hold:
	\begin{itemize}
		\item[\rm(a)] There exists $N\in\NN$ such that for $n\ge N$, $-\frac{1}{2}I+P_n(T_n{\Kcal_{\Om}^\omega})P_n$ has, counting multiplicity, exactly $\mu$ characteristic values in some closed ball $\overline{B_{\delta}(\omega_0)}$ for constant $\delta>0$. 
		
		\item[\rm(b)] We denote by $\tomega_{n,j}$ the characteristic values of $-\frac{1}{2}I+P_n(T_n{\Kcal_{\Om}^\omega})P_n$ for $j=1,\dots,\mu$. For each $j$, the sequence $\{\widetilde{\omega}_{n,j}\}$ converges to $\omega_0$ as $n\to\infty$.
		
		\item[\rm(c)] If $\p\Om$ is of class $C^{3,\alpha}$ for some $1/2<\alpha<1$, there exists a constant $C>0$ such that
		\beq\notag
		\left|\widetilde{\omega}_{n,j}-\omega_0\right|\le C {n^{-\frac{3}{2}}},\quad\mbox{for all }n\ge N,\ j=1,\dots,\mu.
		\eeq
	\end{itemize}
\end{theorem}

It is worth noting that one can numerically compute the Laplacian eigenvalues with this matrix by Theorem \ref{theorem:main:second}, neither requiring the mesh generation nor relying on the  quadrature rules for numerical integration (refer to subsection \ref{subsection:K:matrix} for the explicit matrix expressions for $P_n(T_n{\Kcal_{\Om}^\omega})P_n$).

\smallskip

An interesting problem regarding the Laplacian eigenvalues is to find the extremal points (that is, domains) and the critical points of functions that depend on the Laplacian eigenvalues of planar points \cite{Henrot:2006:EPE}. One may employ numerical techniques, including the gradient descent method, to solve such problems. Problem (ii)--estimation of shape perturbation of the eigenvalues--needs to be addressed before the numerical computations, both on simple eigenvalues and on eigenvalues with \emph{multiplicity}. The latter naturally arise due to symmetries of the domain and pose a special challenge. 
We show that, generically, under a small perturbation on the shape of size $\eps$, the eigenvalues with multiplicity bifurcate into as many eigenvalues counting with multiplicity, and we compute the first-order change, of size $\eps$, corresponding to the shape derivative.

To state the result in detail let us introduce some terminologies.
The Riemann mapping theorem asserts that there uniquely exist a constant $\gamma>0$ called the {\it conformal radius} of $\Om$ and a conformal mapping $\Psi_\Om$ from $\CC\setminus\overline{\gamma\mathbb{D}}$ onto $\CC\setminus\overline{\Om}$ that admits the Laurent series expansion
\beq\label{eq:Psi:laurent}
\Psi_\Om(\xi)=\xi+\sum_{n=0}^\infty a_n \xi^{-n},\quad |\xi|>\gamma,
\eeq
for $a_n \in \mathbb{C}$. We also write $\Psi$ when $\Om$ is clear from the context and refer to, for instance, \cite{Jung:2021:SEL} for the numerical computation of $a_n$ of a given parametrized curve $\p\Om$.
We then perturb each coefficient in the Laurent series individually and denote the resulting domain by $\Om[k,\ep]$ so that 
\beq\label{Psi:deformed}
\Psi_{\Om[k,\varepsilon]}(\xi)=\Psi_\Omega(\xi) + \varepsilon \xi^{-k},\quad|z|>\gamma.
\eeq
We derive asymptotic relations for the layer potential operators depending on the shape deformation of the domain. Again, by applying the Gohberg--Sigal theory to the asymptotic relations, we obtain the following shape derivative formula for the Laplacian eigenvalues. Our main result for problem (ii) is as follows:
\begin{theorem}\label{thm:shapeder:Dirichlet:multipleeig}
	Assume that $\Psi_\Omega$ admits a conformal extension to $\mathbb{C}\backslash\overline{\gamma_2\mathbb{D}}$ for some $0< \gamma_2 < \gamma$.
	Let $\lambda=\omega_0^2>0$ be any Dirichlet (resp. Neumann) eigenvalue of $-\Delta$ on $\Omega$ of multiplicity $\mu$. Fix $j\in\mathbb{N}$. There exists $\delta>0$ such that, for any $\eps\in\CC$ with $0<|\eps|<\delta$ and $j=1,\cdots,\mu$, there exists a Dirichlet (resp. Neumann) eigenvalue $\lambda_{j}(\Om[k,\eps])$ of $-\Delta$ on $\Omega[k,\varepsilon]$ satisfying 
	\begin{align}\label{eq:asympt:Dirichlet:simplified}
		\left|\lambda_j(\Om[k,\varepsilon]) - \lambda - 2\left|\varepsilon\right|\eta_j\,\sqrt{\lambda}\right| \le C |\varepsilon|^{1+\frac{1}{\mu_j}},\quad \eta_j=\eta_j(\omega_0,k,\arg(\varepsilon)),
	\end{align}
	where $\eta_j$ and $\mu_j$ are defined depending on the $|\eps|$-order terms in the asymptotic expansion of 
	$\Kcal_{\Om[k,\eps]}^\omega$ by Definition \ref{eta:1_mu} in Subsection \ref{subsec:shapederiv}, $\mu_j$ is the size of the largest Jordan block with eigenvalue $\eta_j$ of a certain matrix, and $C$ is a positive constant independent both of $|\varepsilon|$ and $\arg(\varepsilon)$. 
\end{theorem}

\section{Preliminaries}\label{sec:preliminary}
\subsection{Operator Generalizations of Complex-Theoretic Results}\label{subsect:GS}

The contents of the section appear in the work of I. C. Gohberg and E. I. Sigal in \cite{Gohberg:1971:OEL} and they serve as generalizations of Rouch\'e's Theorem and the Argument Principle of complex function theory to operator-valued functions. See also \cite{Ammari:2009:LPT} for further explanations.

For Banach spaces $\mathcal{B}$ and $\mathcal{B}'$, we denote by $\mathcal{L}(\mathcal{B},\mathcal{B}')$ the algebra of all bounded linear operators from $\mathcal{B}$ to $\mathcal{B}'$. We will be dealing with functions $F(z)$ defined on a subset of $\mathbb{C}$ that are holomorphic in a (possibly deleted) neighborhood of some $z_0 \in \mathbb{C}$ and with values in $\mathcal{L}(\mathcal{B},\mathcal{B}')$. 
We call $z_0$ a {\it characteristic value} of $F(z)$ if there exists a vector-valued function $f(z)$ called the \emph{root function} of $F(z)$ associated with $z_0$ taking values in $\mathcal{B}$ such that 
\begin{itemize}
	\item[\rm i)] $f(z)$ is holomorphic at $z_0$ and $f(z_0) \neq 0$.
	\item[\rm ii)] $F(z)f(z)$ is holomorphic at $z_0$ and $F(z_0)f(z_0)=0$. 
\end{itemize}
Obviously,
$$F(z)f(z)=(z-z_0)^{m(f)}g(z),$$ 
for some $m(f) \in \mathbb{N}$ and a vector-valued function $g(z)$ with values in $\mathcal{B}'$ that is holomorphic at $z_0$ and such that $g(z_0)\ne0$. The vector $f(z_0)\in\mathcal{B}$ is said to be an {\it eigenvector} of $F(z)$ associated with $z_0$ and the closure of all such eigenvectors is the \emph{kernel} of $F(z)$ associated with $z_0$ and is denoted by $\operatorname{Ker}F(z_0)$. 

For an eigenvector $\varphi$ associated with $z_0$, we set
$$\operatorname{rank}(\varphi)=\max\left\{m(f)\,:\,f\mbox{ is a root function satisfying }f(z_0)=\varphi\right\}.$$
When $\dim \operatorname{Ker} F(z_0)=n<\infty$ and $\operatorname{rank} \varphi < \infty$ for every $\varphi \in \operatorname{Ker} F(z_0)$, the \emph{canonical system of eigenvectors} of $F(z)$ associated with $z_0$ consists of the eigenvectors $\varphi_1,\dots,\varphi_n\in\operatorname{Ker}F(z_0)$ such that $\varphi_1$ is of maximum rank and each $\varphi_j, j=2, \cdots, n,$ is of maximum rank in the subspace $\operatorname{Ker} F(z_0) \setminus \operatorname{span}\left\{\varphi_1,\dots,\varphi_{j-1}\right\}$.
We now define the {\it null multiplicity} of a characteristic value $z_0$ of $F(z)$ as $$N(F(z_0))=\sum_{j=1}^n \operatorname{rank}(\varphi_j).$$ For $z_0$ that is not a characteristic value, we set $N(F(z_0))=0$. 

Suppose further that $F^{-1}(z)$ exists and that $z_0$ is a characteristic value of $F^{-1}(z)$. If $\dim \operatorname{Ker} F^{-1}(z_0)=m<\infty$, then we again consider the canonical system of eigenvectors $\phi_1, \cdots, \phi_m \in \operatorname{Ker} F^{-1}(z)$ and define the \emph{polar multiplicity} of $z_0$ as
$$
N(F^{-1}(z_0))=-\sum_{j=1}^m \operatorname{rank}(\phi_j).
$$
We now define the  {\it multiplicity} of the characteristic value $z_0$ of $F(z)$ as 
$$M(F(z_0))=N(F(z_0))-N(F^{-1}(z_0)).$$
If $z_0$ is a characteristic value of $F(z)$ only, then $M(F(z_0))=N(F(z_0))$. On the other hand, if it is a characteristic value of $F^{-1}(z)$ only, then $M(F(z_0))=-N(F^{-1}(z_0))$. Obviously, if $F(z)$ is holomorphic at $z_0$ and $F(z_0)$ is an invertible operator, then $M(F(z_0))=0$. In this case, $z_0$ is called a \emph{regular point} of $F(z)$.

If $z_0$ is a pole of  $F(z)$ so that $F(z)$ can be expressed as a Laurent series
$$
F(z)=\sum_{j=-p} ^{\infty} (z-z_0)^j T_j,
$$
in a neighborhood of $z_0$, then $z_0$ is called a \emph{normal point} of $F(z)$ if the operators $T_j\in \mathcal{L}(\mathcal{B},\mathcal{B}'), -p \leq j \leq -1$, are of finite rank [$F(z)$ is said to be \emph{finitely meromorphic} at $z_0$], $T_0$ is Fredholm [$F(z)$ is said to be \emph{of Fredholm type} at $z_0$], and all points in some deleted neighborhood of $z_0$ are regular. According to \cite[Lemma 2.1]{Gohberg:1971:OEL}, a normal point of $F(z)$ is also a normal point of $F^{-1}(z)$.

For a simply connected domain $\Omega \subset \CC$ with rectifiable boundary $\partial \Omega$, we say that $F(z)$ is {\it normal with respect to} $\partial \Omega$ if
it is finitely meromorphic and of Fredholm type in $\Omega$, continuous on $\partial \Omega$, and invertible at all points $z_0 \in \overline{\Omega}$ except at a finite number of normal points of $F(z)$ in $\Omega$. If $F(z)$ is normal with respect to $\p \Omega$ and $v_j$, $j=1,\cdots,\sigma$, are the characteristic values and poles of $F(z)$ in $\Omega$, then we set
$$\mathcal{M}(F(z);\p \Omega)=\sum_{j=1}^{\sigma} M(F(v_j)).$$

The normality condition on $F(z)$ with respect to $\partial \Omega$ is an important assumption in the generalizations of the Argument Principle and Rouch{\'e}'s Theorem given below to operator-valued functions. The proof of these results relies on a certain factorization of $F(z)$ at a normal point $z_0$ involving rank one projections when the corresponding Fredholm operator $T_0$ has index
$$
\text{Ind }T_0=\dim \operatorname{Ker} T_0 - \dim \operatorname{Ker} T_0^*=0,
$$
and $M(F(z_0))$ can be easily calculated using this factorization (see \cite[Theorem 3.1 and Theorem 3.2]{Gohberg:1971:OEL} for more details). Let $F(z)$ and $A(z)$ be operator-valued functions taking values in $\mathcal{L}(\mathcal{B}, \mathcal{B}')$.

\begin{theorem}[{\rm Gohberg--Sigal \cite[Theorem 2.2]{Gohberg:1971:OEL}}]
	\label{general:rouche}
	Let $F(z)$ be normal with respect to a rectifiable boundary $\p \Omega$ of a simply connected domain $\Omega \subset \mathbb{C}$. If $A(z)$ is finitely meromorphic in $\Omega$, continuous on $\p \Omega$, and satisfies
	$$\left\| F^{-1}(z)\,A(z) \right\|_{\mathcal{L}(\mathcal{B})}<1, \quad z\in\p \Omega,$$
	then $F(z)+A(z)$ is normal with respect to $\p \Omega$ as well and 
	$$\mathcal{M}(F(z)+A(z);\p \Omega)=\mathcal{M}(F(z);\p \Omega).$$
\end{theorem}
\begin{theorem}[{\rm Gohberg--Sigal \cite[Theorem 4.1]{Gohberg:1971:OEL}}]
	\label{arg:principle}
	Let $F(z)$ be normal with respect to a rectifiable boundary $\p \Omega$ of a simply connected domain $\Omega \subset \mathbb{C}$. Denote by $v_j$, $j=1,\cdots,\sigma$, all the characteristic values and poles of $F(z)$ in $\Omega$. For a scalar function $f$ that is analytic in $\Omega$ and continuous on $\overline{\Omega}$,
	$$\frac{1}{2\pi \mathrm{i}}\operatorname{tr}\int_{\p \Omega} f(z)F'(z)F^{-1}(z)dz=\sum_{j=1}^\sigma M(F(v_j))f(v_j),$$
\end{theorem}
where "tr" stands for the trace of a finite-rank operator.

\subsection{Exterior Conformal Mapping and Geometric Density Basis Functions}\label{subsec:conformal}
Let $\Om \subset \mathbb{C}$ be a simply connected bounded $C^{1,\alpha}$ domain for some $\alpha>0$.
Recall that there is the exterior conformal mapping $\Psi$ associated with $\Om$ given by \eqnref{eq:Psi:laurent}.
As $\p\Om$ is assumed to be of class $C^{1,\alpha}$, $\Psi'$ also continuously extends to $\p\Om$ by the Kellogg--Warschawski Theorem (see \cite[Theorem 3.6]{Pommerenke:1992:BCM}).
Since $\xi$ can be expressed as 
\beq\label{coord:def}
\xi=e^{\rho+\mathrm{i}\theta},
\eeq
for some $(\rho,\theta)\in[\log \gamma,\infty)\times[0,2\pi)$, the length element $d\sigma$ on $\p\Om$ satisfies $d\sigma(z)=h(z)d\theta$ with the scale factor $h(z)= \left|\frac{\p \Psi}{\p \rho}\right|=\left|\frac{\p \Psi}{\p \theta}\right|.$ 
For a function $f$ on $\p\Om$, we have
\beq\label{pOm:normalderi}
\frac{\p f}{\p\nu}\Big|^{+}_{\p\Om}(z)=\frac{1}{h(z)}\lim_{\rho\to(\log \gamma)^+}\frac{\p f(\Psi(e^{\rho+i\theta}))}{\p\rho}\quad\mbox{on }\p\Om.
\eeq

We denote by $F_m(z)$, $m\in\ZZ$, the \emph{Faber polynomials} associated with $\Psi$ \cite{Faber:1903:UPE}. 
A notable property of the Faber polynomials is that $F_m(\Psi(\xi))$ has only one positive order term $\xi^m$, as shown in
\beq\label{Faber:gener}
F_m(\Psi(\xi))=\xi^m+\sum_{k=1}^\infty c_{m,k}\xi^{-k},\quad |\xi|>\gamma,
\eeq
where the constants $c_{m,k}$ are the \emph{Grunsky coefficients}. These coefficients can be easily computed using the recurrence relation (see, for instance,  \cite{Duren:1983:UF, Grunsky:1939:KFS}) $c_{1,m} = a_m$, $c_{m,1} = ma_m$, and
\beq\label{Grunsky:compute:formula}
\begin{aligned}
	&c_{m,n+1} = c_{m+1,n} - a_{m+n} + \sum_{s=1}^{m-1} a_{m-s}c_{s,n} - \sum_{s=1}^{n-1} a_{n-s}c_{m,s},\quad m,n\geq1.
\end{aligned}
\eeq
It can be shown that 
\begin{equation}\label{Grun:inter}
	kc_{m,k} =mc_{k,m}\mbox{ for all }m,k\in\mathbb{N} .
\end{equation}
Further, for $\p\Om$ that is of class $C^{1+p,\alpha}$ for some $p \in \{0\} \cup \mathbb{N}$ and $\alpha\in(0,1)$, there exists a constant $M$ independent of $s$ and $t$ such that (see \cite[Lemma 1.5]{Suetin:1974:POB})
$$
\left|\sum_{k=1}^\infty \frac{c_{s,k}\,\overline{c_{k,t}}}{\gamma^{s+t+2k}}\right|\le\frac{M}{s^{p+\alpha}\,t^{p+\alpha}}\ \mbox{for all }s,t\in\mathbb{N}.
$$
Combining this bound with \eqnref{Grun:inter}, we obtain the following result:
\begin{theorem} \label{Suetin}
	If  $\p\Om$ is of class $C^{3,\alpha}$ for some $\alpha>0$, then
	$$\frac{c_{s,k}\,\overline{c_{k,s}}}{\gamma^{2s+2k}}= \frac{c_{k,s}\,\overline{c_{s,k}}}{\gamma^{2s+2k}}\le \frac{M}{ \left(\max\{s,k\}\right)^{4+2\alpha}}, \quad\mbox{for all }s,k\in\mathbb{N}.
	$$
\end{theorem}

For $z=\Psi(\gamma e^{i\theta})$, we now set $\widetilde{\eta}_0(z)=\eta_0(z)=1$, $\widetilde{\zeta}_0(z)=\zeta_0(z)=\frac{1}{h(z)}$, and
\begin{align*}
	\widetilde{\eta}_m(z)&=e^{\mathrm{i}m\theta},\quad \widetilde{\zeta}_m(z)=\frac{e^{im\theta}}{h(z)},\\
	\eta_{m}(z)&=\frac{1}{\sqrt{|m|}}\,e^{ \mathrm{i} m\theta},\quad
	\zeta_{m}(z)=\sqrt{|m|}\,\frac{e^{\mathrm{i}m\theta}}{h(z)}\quad\mbox{for }m\in\ZZ\setminus\{0\}.
\end{align*}
Since $\p\Om$ is of class $C^{1,\alpha}$ for some $\alpha>0$, both $h$ and $1/h$ are non-vanishing and continuous on $\p\Om$.
The work by Y. Jung and the fourth author shows that the set $\{\widetilde{\eta}_{m}\}_{m\in\ZZ}$ forms a basis for $L^{2}(\p\Omega)$ (\cite{Jung:2021:SEL}) so that for $\varphi\in L^{2}(\p\Om)$,
$$\label{newexp:coeff}
\varphi=\sum_{m\in\ZZ}a_m\widetilde{\eta}_m,\quad\mbox{where }a_m=\frac{1}{2\pi}\int_{\p\Om}\varphi\,\widetilde{\zeta}_{-m}\, d\sigma.
$$
An inner product on $L^{2}(\p\Om)$ is given in the usual way as
\beq
\label{def:newinner} 
\Big\langle\sum_{m\in\mathbb{Z}}\alpha_m\widetilde{\eta}_m,\, \sum_{m\in\mathbb{Z}}\beta_m\widetilde{\eta}_m\Big\rangle_{0}:=\sum_{m\in\mathbb{Z}}{\alpha_m}\overline{\beta_m},
\eeq
and the norm $\|\cdot\|_{0}$ associated with $\langle\cdot,\cdot\rangle_{0}$ is equivalent to the standard norm on $L^{2}(\p\Om)$.  

The layer potentials of the Laplacian admit the following expansions in terms of these geometric basis functions:
\begin{theorem}[\cite{Jung:2021:SEL}]\label{lem:uss}
	On $\p\Om$, $\mathcal{S}_\Omega^0[\zeta_0]=\log \gamma$, $\mathcal{K}_\Omega^0[\eta_0]=\frac{1}{2}$, 
	and
	\begin{align}\notag
		\mathcal{S}_\Omega^0\left[\zeta_m\right]&=-\frac{1}{2}\left(\eta_m+\sum_{k=1}^\infty \frac{\sqrt{k}}{\sqrt{m}}\frac{c_{m,k}}{\gamma^{m+k}}\,\eta_{-k}\right),\quad
		\mathcal{S}_\Omega^0\left[\zeta_{-m}\right]=-\frac{1}{2}\left(\eta_{-m}+\sum_{k=1}^\infty\frac{\sqrt{k}}{\sqrt{m}}\frac{\overline{c_{m,k}}}{\gamma^{m+k}}\,\eta_{k}\right),\\\notag
		\mathcal{K}_\Omega^0\left[\eta_m\right]&=\frac{1}{2}\sum_{k=1}^\infty \frac{\sqrt{k}}{\sqrt{m}}\frac{c_{m,k}}{\gamma^{m+k}}\,\eta_{-k},\quad
		\mathcal{K}_\Omega^0\left[\eta_{-m}\right]=\frac{1}{2}\sum_{k=1}^\infty\frac{\sqrt{k}}{\sqrt{m}} \frac{\overline{c_{m,k}}}{\gamma^{m+k}}\,\eta_{k}\quad\mbox{for }m \in \mathbb{N}.
	\end{align}
\end{theorem}

\section{Approximations of $\Kcal_\Om^\omega$ by Finite Matrices}
We approximate the Laplacian eigenvalues by the eigenvalues of certain finite-dimensional matrices that are ''close" to the operator $\Kcal_\Om^\omega$ and whose entries are explicitly determined by the conformal mapping coefficients of $\Om$. Let $\Om$ be a $C^{1, \alpha}$ domain for some $\alpha>0$.

\subsection{Finite Section Method for $\Kcal_\Om^\omega$}
In view of (\ref{eq:fundamental:sol}) and (\ref{Kdef}), we express $\Kcal_\Om^\omega$ as
$$
\ds\Kcal_{\Om}^\omega=\lim_{n\rightarrow\infty} T_n\Kcal_\Om^\omega,
$$
where
\beq\label{Kcal_k_exp}
\begin{aligned}
	&\ds T_n\Kcal_{\Om}^\omega = \Kcal_\Om^0 + \sum_{j=1}^n \left(\omega^{2j}\log \omega\right) \Kcal_{\Om}^{(j,1)}
	+\sum_{j=1}^n \omega^{2j} \Kcal_{\Om}^{(j,2)},\\
	&\ds\mathcal{K}_{\Omega}^{(j,1)}[\varphi](x):=\int_{\p\Omega}b_j\frac{\p |x-y|^{2j}}{\p\nu_y}\varphi(y)\,d\sigma(y),\\
	&\ds\mathcal{K}_{\Omega}^{(j,2)}[\varphi](x):=\int_{\p\Omega}\frac{\p \left(|x-y|^{2j} b_j(\log|x-y|+c_j)\right)}{\p\nu_y}\varphi(y)\,d\sigma(y),
\end{aligned}
\eeq
and the constants $b_j$ and $c_j$ are given by (\ref{bjcj}).
For notational convenience, we set
$\Kcal_{\Om}^{(0,1)}=0$ and $\Kcal_{\Om}^{(0,2)}=\Kcal_{\Om}^0$. The following lemmas are essential in deriving error estimates for the computation of the Laplacian eigenvalues:
\begin{lemma}\label{lemma:kminustk:pointwise}
	For a fixed $\omega\in\CC\setminus \RR^{-}$, 
	$$\lim_{n\rightarrow\infty} \|\mathcal{K}_\Omega^\omega-T_n\mathcal{K}_{\Omega}^\omega\big\|_{\mathcal{L}(L^2(\p\Om))}=0.$$
\end{lemma}
\begin{proof}
	We first note that 
	\beq\label{Kcal_k_exp:2}
	\begin{aligned}
		&\ds \left | \frac{\p |x-y|^{2j}}{\p\nu_y} \right|=2j|\nu_y \cdot (y-x)||x-y|^{2j-2} \leq 2j \max_{x \in \partial \Om} (1+2|x|)^{2j-1} ,\\
		&\ds \left| \frac{\p \left(|x-y|^{2j} \log|x-y| \right)}{\partial \nu_y} \right| \leq 2j(|x-y|^{2j-2}+|x-y|^{2j})+|x-y|^{2j-1} \leq 5j  \max_{x \in \partial \Om} (1+2|x|)^{2j}.
	\end{aligned}
	\eeq
	Setting $K=\max_{x \in \partial \Om} (1+2|x|)$, we then have
	\beq \label{Kcalbounds} 
	\left \Vert {\mathcal{K}_{\Omega}^{(j,1)}[\varphi] }\right \Vert_{L^2(\partial \Omega)}, \left \Vert{\mathcal{K}_{\Omega}^{(j,2)}[\varphi] }\right \Vert_{L^2(\partial \Omega)} \leq 5j |b_j|(1+|c_j|)K^{2j}|\partial \Om|^{\frac{1}{2}}\|\varphi\|_{L^2(\partial \Omega)},
	\eeq
	and this, in turn, gives
	\begin{equation}\label{inequality:KmTK}
		\left \Vert \mathcal{K}_{\Omega}^{\omega}-T_n\mathcal{K}_{\Omega}^{\omega} \right \Vert_{\mathcal{L}(L^2(\partial \Omega))} \leq C \sum_{j=n+1}^{\infty} |\omega|^{2j} (1+\log |\omega|)\frac{\log j}{2^{2j}(j!)^2}K^{2j},
	\end{equation}
	for some constant $C$.
	We observe that the right-hand side of \eqref{inequality:KmTK} converges to zero as $n\to\infty$.
\end{proof}

In what follows, let $P_n$ denote the orthogonal projection from $L^2(\p\Om)$ onto $\mbox{span }\{\widetilde{\eta}_m:\, -n\leq m\leq n\}$. 

\begin{lemma}\label{lemma:finitesection:conv}
	For a fixed $\omega\in\CC\setminus \RR^{-}$, 
	$$\lim_{n\rightarrow\infty} \|\mathcal{K}_\Omega^\omega-P_n\mathcal{K}_{\Omega}^\omega{P_n}\big\|_{\mathcal{L}(L^2(\p\Om))}=0.$$
	
\end{lemma}
\begin{proof}
	Since $\partial \Omega$ is of class $C^{1, \alpha}$, both $\mathcal{K}_{\Omega}^{\omega}$ and $\mathcal{K}_{\Omega}^{{\omega}*}$ are compact operators on $L^2(\partial \Omega)$. Let 
	$$
	\mathcal{B}=\{f \in L^2(\partial \Omega): \|f\|_0=1\},
	$$
	using the norm $\|\cdot \|_0$ associated with (\ref{def:newinner}),
	and define a sequence of continuous functions $r_n: \overline{\mathcal{K}_\Omega^{\omega}(\mathcal{B})} \rightarrow \mathbb{R}$ by
	$$
	r_n(f)= \| (I-P_n)f\|_0.
	$$
	For each fixed $f \in \overline{\mathcal{K}_\Omega^{\omega}(\mathcal{B})}$, $r_n(f)$ is decreasing in $n$ and $\lim_{n \rightarrow \infty} r_n(f)=0.$ Since $r_n(f) \rightarrow 0$ uniformly on $\overline{\mathcal{K}_\Omega^{\omega}(\mathcal{B})}$, 
	$$
	\lim_{n \rightarrow \infty} \|(I-P_n)\mathcal{K}_\Omega^{\omega}[f] \|_0=0 \text{ uniformly for }f \in \mathcal{B}.
	$$
	Similarly,
	$$
	\lim_{n \rightarrow \infty} \|(I-P_n)(\mathcal{K}_\Omega^{\omega})^*[f] \|_0=0 \text{ uniformly for }f \in \mathcal{B}.
	$$
	
	Denoting by $\| \cdot \|_0$ the operator norm with respect to the norm $\| \cdot \|_0$ of $L^2(\partial \Omega)$, we then have
	$$
	\begin{aligned}
		\|\mathcal{K}_\Omega^\omega-P_n\mathcal{K}_{\Omega}^\omega{P_n}\big\|_0 
		& \leq \|(I-P_n)\mathcal{K}_{\Omega}^{\omega}\|_0+\|P_n \mathcal{K}_{\Omega}^{\omega}(I-P_n)\|_0 \\
		& \leq \|(I-P_n)\mathcal{K}_{\Omega}^{\omega}\|_0+\|(I-P_n)\mathcal{K}_{\Omega}^{\omega*}\|_0 \\
		& \rightarrow 0.
	\end{aligned}
	$$
	Finally, we note that the norm $\| \cdot \|_0$ is equivalent to the standard norm on $L^2(\partial \Omega)$.
\end{proof}
\begin{lemma} \label{Kconvergence}
	Let $B_{\delta}$ be a disk of radius $\delta$ that is compactly contained in $\CC\setminus\sqrt{-1}\RR^-$. Then
	\beq\notag
	\sup_{\omega\in\partial B_{\delta}}
	\big\|\mathcal{K}_\Omega^\omega-P_n (T_n{\Kcal_{\Om}^\omega}) {P}_n\big\|_{\mathcal{L}(L^2(\p\Om))}
	\to0\quad\mbox{as }n\to\infty.
	\eeq
\end{lemma}
\begin{proof}
	For $n \in \mathbb{N}$ and $w \in \partial B_{\delta}$, we let
	$$
	f_n(w):=\big\|\mathcal{K}_\Omega^\omega-P_n (T_n{\Kcal_{\Om}^\omega}) {P}_n\big\|_{\mathcal{L}(L^2(\p\Om))}.
	$$
	Applying the previous lemmas to 
	$$
	\mathcal{K}_\Omega^\omega-P_n (T_n{\Kcal_{\Om}^\omega}) {P}_n=\mathcal{K}_\Omega^\omega-P_n\Kcal_{\Om}^\omega
	P_n+P_n(\Kcal_{\Om}^\omega-T_n\Kcal_{\Om}^\omega)P_n,
	$$
	we obtain for a fixed $w$,
	$$
	\lim_{n \rightarrow \infty} \|\mathcal{K}_\Omega^\omega-P_n (T_n{\Kcal_{\Om}^\omega}) {P}_n\|=0.
	$$
	
	In the remainder of the proof, we show that the convergence is, in fact, uniform for $w \in \partial B_{\delta}$. The equivalence of the norms $\| \cdot \|_{\mathcal{L}(L^2(\partial \Omega))}$ and $\| \cdot \|_0$ first yields the bound
	$$
	\|P_n\|_{\mathcal{L}(L^2(\partial \Omega))} \leq C_1\|P_n\|_0=C_1,
	$$
	for some constant $C_1$ independent of $n$. Then for $\omega, \widetilde{\omega} \in \partial B_{\delta}$,
	\begin{align*}
		\left|f_n(\omega)-f_n(\widetilde{\omega})\right|
		&\le \big\|\mathcal{K}_{\Omega}^\omega-\mathcal{K}_{\Omega}^{\widetilde{\omega}}\big\|
		+\big\|P_n\big({T_n\Kcal_{\Om}^\omega} -T_n\mathcal{K}_{\Omega}^{\widetilde{\omega}}\big) P_n\big\|
		\le \big\|\mathcal{K}_{\Omega}^\omega-\mathcal{K}_{\Omega}^{\widetilde{\omega}}\big\|+C_1^2\big\|{T_n\Kcal_{\Om}^\omega}-T_n\mathcal{K}_{\Omega}^{\widetilde{\omega}} \big\|\\
		&\le (1+C_1^2)\sum_{j=1}^\infty\left(\left|\omega^{2j}-\widetilde{\omega}^{2j}\right|\,\big\|\Kcal_{\Om}^{(j,2)}\big\|+\left|\omega^{2j}\log\omega-\widetilde{\omega}^{2j}\log\widetilde{\omega}\right|\,\big\|\Kcal_{\Om}^{(j,1)}\big\|\right)\\
		&\leq (1+C_1^2)\left(C_2\left|\omega-\widetilde{w}\right|+C_3\left|\log\omega-\log\widetilde{w}\right|\right),
	\end{align*}
	with 
	\beq
	\begin{aligned}
		\label{Lemma3.3C2}
		C_2&=\sum_{j=1}^\infty\, 2j\Big(1+\sup_{z\in\partial B_\delta}|z|\Big)^{2j-1}\Big(\big\|\Kcal_{\Om}^{(j,2)}\big\|+\sup_{z\in\partial B_\delta}\left|\log z\right|\,\big\|\Kcal_{\Om}^{(j,1)}\big\|\Big)<\infty,\\
		C_3&=\sum_{j=1}^\infty\sup_{z\in\partial B_\delta}|z|^{2j}\, \|\Kcal_{\Om}^{(j,1)}\|<\infty.
	\end{aligned}
	\eeq
	Note that the convergence of the series in \eqnref{Lemma3.3C2} is justified by \eqnref{Kcalbounds}. Now, since $\p B_\delta$ is bounded away from the origin, and $\left|\log\omega-\log\widetilde{\omega}\right|/|\omega-\widetilde{\omega}|$ is bounded uniformly for $\omega,\,\widetilde{\omega}\in\partial B_\delta$, 
	we have
	$$\left|f_n(\omega)-f_n(\widetilde{\omega})\right|\leq C_4 \left|\omega-\widetilde{w}\right|,$$
	for a constant $C_4$ independent of $n$ and $\omega,\widetilde{\omega}\in \p B_\delta$. 
	Hence, $\{f_n\}$ is equicontinuous on the compact set $\p B_\delta$. Using the fact that $f_n$ converges to $0$, we conclude that $f_n$ uniformly converges to $0$ on $\p B_\delta$. 
\end{proof}

\subsection{Convergence Estimate assuming a greater regularity of $\p\Om$}
We derive the convergence rate for the operator in Lemma \ref{lemma:finitesection:conv}, assuming a greater regularity than $C^{1,\alpha}$ of $\p\Om$, as follows:
\begin{theorem}\label{theorem:main:smooth}
	Let $\p\Om$ be of class $C^{3,\alpha}$ for some $1/2<\alpha<1$. 
	We have
	$$
	\big\|\mathcal{K}_\Om^\omega - P_n \mathcal{K}_\Om^\omega P_n\big\|_{\mathcal{L}(L^2(\p\Om))}\le\frac{R(\omega)}{n^{3/2}},\quad\mbox{for all }\omega\in\mathbb{C}\backslash\sqrt{-1}\mathbb{R}^-\mbox{ and } n\in\NN,
	$$
	where $R(\omega)$ is a continuous function for $\omega\in\mathbb{C}\backslash\sqrt{-1}\mathbb{R}^-$.
\end{theorem}

It suffices to prove the theorem with the operator norm with respect to the norm $\| \cdot \|_0$ of $L^2(\partial \Omega)$. Moreover, one can decompose $\mathcal{K}_{\Omega}^{\omega}$ on $L^2(\partial \Omega)$ into
\beq
\label{Kcaldecomp}
\Kcal_{\Om}^\omega =\Kcal_{\Om}^0 + \mathcal{R}_1+\mathcal{R}_2, \text{ where }
\eeq
\begin{align}\label{R1R2}
	\mathcal{R}_1[\varphi](z)
	&=-\int_0^{2\pi} \bigg( \sum_{j=1}^\infty \omega^{2j}\Big((\log\omega)jb_j+jb_j c_j+\frac{1}{2}b_j\Big) 
	\big|z-\Psi(\gamma e^{\mathrm{i}\theta})\big|^{2j-2}\, p(\theta;z)\bigg)\varphi\big(\Psi(\gamma e^{\mathrm{i}\theta})
	\big)\,d\theta,\\ 
	\mathcal{R}_2[\varphi](z)
	&=-\int_0^{2\pi} \log \left|z-\Psi(\gamma e^{\mathrm{i}\theta})\right|
	\bigg(\sum_{j=1}^\infty \omega^{2j} jb_j 
	\big|z-\Psi(\gamma e^{\mathrm{i}\theta})\big|^{2j-2}\, p(\theta;z)\bigg)\varphi\big(\Psi(\gamma e^{\mathrm{i}\theta})
	\big)\,d\theta, \notag
\end{align}
and 
\beq\notag
p(\theta;z)=\Psi'(\xi)\xi\,\overline{\left(z-\Psi(\xi)\right)}+\overline{\Psi'(\xi)\,\xi}\left(z-\Psi(\xi)\right),\quad \xi=\gamma e^{\mathrm{i}\theta}.
\eeq

This can be easily seen by noting that for $z\in\p\Om$ and $\tau=\Psi(\xi)$ with $\xi=e^{\rho+\mathrm{i}\theta}$,
we have from \eqnref{pOm:normalderi} that 
\begin{align*} 
	& h(\tau)\,\frac{\p (|z-\tau|^{2j})}{\p\nu_\tau}\bigg|_{\p\Om}^{+}(\tau)
	=\lim_{\rho\to(\log \gamma)+} \frac{\p (z-\tau)^j(\overline{z-\tau})^j}{\p \rho}
	=-j\Big|z-\Psi(\gamma e^{\mathrm{i}\theta})\Big|^{2j-2} p(\theta;z),\\
	&h(\tau)\,\frac{\p \log |z-\tau|}{\p\nu_\tau}\bigg|_{\p\Om}^{+}(\tau)
	=\lim_{\rho\to(\log \gamma)+}\frac{1}{2}\frac{\p \ln\big((z-\tau)(\overline{z-\tau})\big)}{\p\rho}
	= -\frac{1}{2}\Big|z-\Psi(\gamma e^{\mathrm{i}\theta})\Big|^{-2} p(\theta;z).
\end{align*}
Using these calculations in the definitions given in  \eqnref{Kcal_k_exp}, we then obtain
\begin{align*}
	&\Kcal_\Om^{(j,1)}[\varphi](z)= - \int_0^{2\pi} j\, b_j 
	\big|z-\Psi(\gamma e^{\mathrm{i}\theta})\big|^{2j-2}\, p(\theta;z)\,\varphi\big(\Psi(\gamma e^{\mathrm{i}\theta})
	\big)d\theta,\\
	&\Kcal_\Om^{(j,2)}[\varphi](z)
	=-\int_0^{2\pi} \big(j\,b_j \log \big|z-\Psi(\gamma e^{\mathrm{i}\theta})\big| + j\,b_j c_j+\frac{1}{2}b_j\big) \big|z-\Psi(\gamma e^{\mathrm{i}\theta})\big|^{2j-2}\, p(\theta;z)\,\varphi\big(\Psi(\gamma e^{\mathrm{i}\theta})
	\big)d\theta.
\end{align*}

\subsubsection{Relevant Lemmas}

We will present a series of lemmas that will be used to estimate the terms in (\ref{Kcaldecomp}). For notational simplicity, we let
$$\label{kstar:def}
n_*=\max({|n|},1), \quad\mbox{for }n\in\ZZ.
$$

\begin{lemma} \label{firstlemma}
	\notag
	Let $x,y>0$ with $\max\{x,y\}>1$. There exists a constant $C_{x,y}$ such that
	\begin{align}
		\sum\nolimits_{q\in\mathbb{Z}}q_*^{-x}\, (q-m)_*^{-y}&\le C_{x,y}\,{m_*^{-\min\{x,y\}}},\quad\mbox{for all }m\in\mathbb{Z}.
	\end{align}
\end{lemma}
\begin{proof}
	Since both sides of the inequality are invariant under the transformation $m\mapsto -m$, we may assume that $m\geq 0$. 
	The result obviously holds for $m=0$, and therefore, by considering the transformation $q\rightarrow -q+m$, it suffices to prove it under the assumption that 
	$m>0$ and $y\geq x$ (and thus, $y>1$). We split the sum into
	\beq\notag
	\sum_{q\in\mathbb{Z}}\frac{1}{q_*^x\, (q-m)_*^y}= \frac{1}{m^y}+\frac{1}{m^x}+\bigg(\sum_{q=-\infty}^{-1}+\sum_{q=1}^{m-1}+\sum_{q=m+1}^\infty\bigg)\frac{1}{q_*^x\, (q-m)_*^y}
	=:\frac{1}{m^y}+\frac{1}{m^x}+I+II+III,
	\eeq
	where the first two terms correspond to $q=0$ and $q=m$, respectively.
	The estimates
	$
	II\leq (\frac{2}{m})^{x}\frac{2y}{y-1}\mbox{ and }
	I, III \leq (\frac{1}{m})^{x}\frac{y}{y-1}
	$
	can then be easily obtained. 
\end{proof}


\begin{lemma}\label{R1}
	For $\varphi=\sum_{n\in\ZZ}\langle \varphi,\widetilde{\eta}_n\rangle_{0}\, \widetilde{\eta}_n \in L^2(\p\Om),$
	\begin{align*}
		\mathcal{R}_1[\varphi]
		&=\sum\nolimits_{m\in\ZZ}\sum\nolimits_{n\in\ZZ}
		\langle \varphi,\widetilde{\eta}_n\rangle_{0}
		\langle\mathcal{R}_1[\widetilde{\eta}_n], \widetilde{\eta}_m\rangle_{0}\, \widetilde{\eta}_m,
	\end{align*}
	$$\label{eq:lemma:R1coeff}
	\big|\langle\mathcal{R}_1[\widetilde{\eta}_n], \widetilde{\eta}_m\rangle_{0}\big|
	\leq  C(\omega){(1+|m|)^{-3}}{(1+|n|)^{-2}},\quad\mbox{for }m,n\in\ZZ,
	$$
	where $C(\omega)$ depends continuously on $\omega\in \CC\setminus\sqrt{-1}\RR^{-}$.
\end{lemma}
\begin{proof}
	
	For $m,n\in\ZZ$, 
	\begin{align*}
		\langle\mathcal{R}_1[\widetilde{\eta}_n], \widetilde{\eta}_m\rangle_{0}
		&=\frac{1}{2\pi} \int_{\p\Om}\mathcal{R}_1[\widetilde{\eta}_n](z)\,\widetilde{\zeta}_{-m}(z)\,d\sigma(z)
		=\frac{1}{2\pi}\int_0^{2\pi}\mathcal{R}_1[\widetilde{\eta}_n](\Psi(\gamma e^{\mathrm{i}\theta_1}))\, e^{-\mathrm{i}m\theta_1}\,d\theta_1.
	\end{align*}
	Substituting $\varphi=\widetilde{\eta}_n$ in the definition of $\mathcal{R}_1$ given by \eqnref{R1R2}, we have
	\begin{align}\label{fourier:inner} 
		\langle\mathcal{R}_1[\widetilde{\eta}_n], \widetilde{\eta}_m\rangle_{0}
		&=\frac{-1}{2\pi}\int_0^{2\pi} \int_0^{2\pi}f(\theta_1,\theta_2;\omega)
		\, e^{\mathrm{i}n\theta_2}\, e^{-\mathrm{i}m\theta_1}\,d\theta_2\,d\theta_1, \text{ where}\\
		\label{def:f} \notag
		f(\theta_1,\theta_2;\omega)&=
		\sum_{j=1}^\infty \omega^{2j}\Big((\log\omega)jb_j+jb_j c_j+\frac{1}{2}b_j\Big) 
		\big|\Psi(\gamma e^{\mathrm{i}\theta_1})-\Psi(\gamma e^{\mathrm{i}\theta_2})\big|^{2j-2}\, p(\theta_2;\Psi(\gamma e^{\mathrm{i}\theta_1})),
	\end{align}
	and in the remainder of the proof, we aim to estimate this term.
	
	First, by the definition of $p(\theta_2;z)$, 
	\beq\notag
	\big|\Psi(\gamma e^{\mathrm{i}\theta_1})-\Psi(\gamma e^{\mathrm{i}\theta_2})\big|^{2j-2}\, p(\theta_2;\Psi(\gamma e^{\mathrm{i}\theta_1}))=\frac{{\rm i}}{j}\,\partial_{\theta_2}\left(\big|\Psi(\gamma e^{\mathrm{i}\theta_1})-\Psi(\gamma e^{\mathrm{i}\theta_2})\big|^{2j}\right).
	\eeq
	Moreover, a term-by-term differentiation of $f(\theta_1, \theta_2; \omega)$ results in
	$$
	\p_{\theta_1}^{3}\,\p_{\theta_2}^{2}\,f
	= \sum_{j=1}^\infty \omega^{2j}\Big((\log\omega)jb_j+jb_j c_j+\frac{1}{2}b_j\Big) 
	\frac{{\rm i}}{j}\,\partial_{\theta_1}^3\,\partial_{\theta_2}^3\left(\big|\Psi(\gamma e^{\mathrm{i}\theta_1})-\Psi(\gamma e^{\mathrm{i}\theta_2})\big|^{2j}\right), 
	$$ 
	so that 
	$$
	\left|\p_{\theta_1}^{3}\,\p_{\theta_2}^{2}\,f\right|
	\leq \sum_{j=1}^\infty |\omega|^{2j}\Big|(\log\omega)jb_j+j b_j c_j+\frac{1}{2}b_j\Big|\,\frac{2^{2j}}{j}\, (2j)^{6}\widetilde{M}^{2j}=:C(\omega),
	$$
	where $\widetilde{M}:=1+\max_{0\le k\le 3}\max_{\theta\in T}\left|\p_{\theta}^{k}\,{\Psi(\gamma e^{{\rm i}\theta})}\right|<\infty$. This is due to the fact that the regularity assumption on $\Psi$ yields
	$$\label{regularity:f}
	\p_{\theta}^{k}\,\Psi(\gamma e^{{\rm i}\theta}),\ \p_{\theta}^{k}\,\overline{\Psi(\gamma e^{{\rm i}\theta})}\in C^{0,\alpha}(T)\quad\mbox{for any } 0\le k\le 3,
	$$
	where $T$ is the one-dimensional torus $[0,2\pi].$

	Since \eqnref{bjcj} shows that $b_j$ and $b_j c_j$ decay fast enough, $C(\omega)$ is a continuous function of $\omega\in\mathbb{C}\backslash\sqrt{-1}\mathbb{R}^-$.
	Moreover, 
	$$
	|\widehat{(\p_{\theta_1}^{3}\,\p_{\theta_2}^{2}\,f)}_{m,n}|=\left|{(2\pi)^{-2}}\int_0^{2\pi}\int_0^{2\pi} (\p_{\theta_1}^{3}\,\p_{\theta_2}^{2}\,f) e^{-{\rm i}m\theta_2} e^{-{\rm i}n\theta_1}\,d\theta_2\,d\theta_1\right|
	\le C(\omega),\quad m,n\in\mathbb{Z}.
	$$
	Finally, \eqnref{fourier:inner} and the fact that $|\widehat{(\p_{\theta_1}^{3}\,\p_{\theta_2}^{2}\,f)}_{m,n}|=|m|^{3}|n|^{2}|\widehat{f}_{m,n}|$ yield the desired result.
\end{proof}

We next have
\begin{align}\notag
	\mathcal{R}_2[\varphi](\Psi(\gamma e^{\rm i \theta_1}))
	&=-\int_0^{2\pi} \log \big|\Psi(\gamma e^{\rm i \theta_1})-\Psi(\gamma e^{\mathrm{i}\theta_2})\big|\,
	g(\theta_1,\theta_2;\omega)\,\varphi\big(\Psi(\gamma e^{\mathrm{i}\theta_2})
	\big)\,d\theta_2,  \text{ where}\\ \label{eq2:part2to1} \notag
	g(\theta_1,\theta_2;\omega)&=\sum_{j=1}^\infty \omega^{2j} j\,b_j 
	\big|\Psi(\gamma e^{\mathrm{i}\theta_1})-\Psi(\gamma e^{\mathrm{i}\theta_2})\big|^{2j-2}\, p(\theta_2;\Psi(\gamma e^{\mathrm{i}\theta_1})).
\end{align}
As in the case for $\mathcal{R}_1$,
\beq\label{eq:fouriera:decay} \notag
g(\theta_1,\theta_2;\omega)=\sum\nolimits_{s,t\in\ZZ}\widehat{g}_{s,t}(\omega)\,e^{{\rm i}s\theta_1}e^{{\rm i}t\theta_2},\quad \left|\widehat{g}_{s,t}\right| \leq C(\omega)\, s_*^{-3} \,{t_*^{-2}} ,
\eeq
for a continuous function $C(\omega)$ independent of $s$ and $t$. 
We then have
\begin{align*}
	g(\theta_1,\theta_2;\omega)\,\widetilde{\eta}_n(\Psi(\gamma e^{{\rm i}\theta_2}))
	=\sum\nolimits_{s,t\in\ZZ}\widehat{g}_{s,t}(\omega)\,e^{{\rm i}{s\theta_1}}e^{{\rm i}(t+n)\theta_2}.
\end{align*}

Since the integral representation of $\mathcal{R}_2$ involves a logarithmic function, we consider the single-layer potential with zero frequency.
We define the constants $S_{s,l}$ by
\beq\label{Snewsum} \notag
\mathcal{S}_\Om^0\left[\widetilde{\zeta}_{s}\right](z)=\sum\nolimits_{l\in\mathbb{Z}} S_{s,l}\,\widetilde{\eta}_l(z).
\eeq
Note that by Theorems \ref{Suetin} and \ref{lem:uss},
\beq\label{S} 
\begin{aligned}
	&S_{0,0}=\log \gamma,\quad S_{s,0}=S_{0,l}=0,\\
	&S_{s,-l}=S_{l,-s}=\overline{S_{-l,s}},\quad\mbox{for } s,l\in\mathbb{N}\backslash\{0\},\\
	&S_{s,l}=-\frac{1}{2|s|}\delta_{s,l},\quad\mbox{for }sl>0,
\end{aligned}
\eeq
and 
\beq\label{S2}
\begin{aligned}
	|S_{s,l}|^2 = \frac{1}{4|s||l|}\frac{c_{|s|,|l|}\,\overline{c_{|l|,|s|}}}{\gamma^{2|s|+2|l|}} \le \frac{M}{4} \frac{1}{|s||l|\left(\max\{|s|,|l|\}\right)^{4+2\alpha}},\quad\mbox{for }sl<0.
\end{aligned}
\eeq

Let $\varphi=\sum_{n\in\ZZ}\beta_n \widetilde{\eta}_n$ with $\sum_{n\in\ZZ}|\beta_n|^2=1$.  
In what follows, $C=C(\omega,\Omega)$ will denote a constant independent of $\beta_n$. Since
\begin{align*}
	\mathcal{R}_2\left[\widetilde{\eta}_n\right](\Psi(\gamma e^{\rm i \theta_1}))
	&= {-2\pi}\,\sum_{s,t\in\ZZ}\widehat{g}_{s,t}\,e^{{\rm i}{s\theta_1}}\,\Scal_{\Om}^0\left[\frac{1}{h(\Psi(\gamma e^{{\rm i}\theta_2}))}e^{{\rm i}(t+n)\theta_2}\right](\Psi(\gamma e^{\rm i \theta_1}))\\
	&={-2\pi}\,\sum_{s,t,l\in\ZZ}\widehat{g}_{s,t}\,S_{t+n,l}\, \widetilde{\eta}_{s+l}(\Psi(\gamma e^{\rm i \theta_1}))\\
	&={-2\pi}\,\sum_{m\in\ZZ}\, \sum_{s,t\in\ZZ}\widehat{g}_{s,t}\,S_{t+n,m-s}\, \widetilde{\eta}_{m}(\Psi(\gamma e^{\rm i \theta_1})),
\end{align*}
it follows that
\begin{align}\notag
	&\frac{1}{4\pi^2}\big\|(\mathcal{R}_2-P_N\mathcal{R}_2P_N)[\varphi]\big\|_{0}^2\\
	\label{eq:part2d1:last:1}
	=&\sum_{|m|>N}\bigg|\sum_{n\in\ZZ} \beta_n \sum_{s,t\in\mathbb{Z}}\widehat{g}_{s,t}\, S_{t+n,m-s}\bigg|^2 + \sum_{|m|\leq N} \bigg|\sum_{|n|>N} \beta_n \sum_{s,t\in\mathbb{Z}}\widehat{g}_{s,t}\,S_{t+n,m-s}\bigg|^2.
\end{align}

We now split the summation $\sum_{s,t\in\mathbb{Z}}\widehat{g}_{s,t}\, S_{t+n,m-s}$ into four cases (i) $t+n=m-s\ne 0$ (ii) $t+n=m-s=0$ (iii) $(t+n)(m-s)<0$ (iv) the rest for which $S_{t+n,m-s}=0$, so that
\begin{align}\notag
	\sum_{s,t\in\mathbb{Z}}\widehat{g}_{s,t}\, S_{t+n,m-s}
	&=\bigg(\sum_{\substack{t+n=m-s\neq 0}} +\sum_{\substack{t+n=m-s=0} } +\sum_{\substack{(t+n)(m-s)<0}} \Bigg)\,\widehat{g}_{s,t}\, S_{t+n,m-s} =: I_1 + I_2 + I_3,
\end{align}
and put
\beq\label{def:Jk12} \notag
\begin{aligned}
	J_k^{(1)}&:=\sum_{|m|>N}\bigg|\sum_{n\in\ZZ} \beta_n \,I_k(m,n)\bigg|^2, \quad
	J_k^{(2)}:=\sum_{|m|\leq N}\bigg|\sum_{|n|>N}\beta_n \,I_k(m,n)\bigg|^2,\quad\mbox{for }k=1,2,3.
\end{aligned}
\eeq
Then
\beq\label{Js}
\frac{1}{4\pi^2}\big\|(\mathcal{R}_2-P_N\mathcal{R}_2P_N)[\varphi]\big\|_{0}^2
\leq 3 \left(J_{1}^{(1)}+J_2^{(1)}+J_3^{(1)}\right)+3 \left(J_{1}^{(2)}+J_2^{(2)}+J_3^{(2)}\right).
\eeq

\begin{lemma} \label{R2estimate}
	\beq \notag
	\frac{1}{4\pi^2}\big\|\mathcal{R}_2-P_N\mathcal{R}_2P_N\big\|_{0}^2=O(N^{-3}).
	\eeq
\end{lemma}
\begin{proof}
	\noindent{\underline{\bf{Estimates for $J_1^{(1)}$} and $J_1^{(2)}$}}.
	Let $t+n=m-s=q \in \ZZ \setminus\{0\}$. First, we see from \eqnref{S} that  $S_{q,q}=-1/(2|q|)$. By substituting $\theta_2=\theta_1$ into the function $g(\theta_1, \theta_2; \omega)$, we have for each $k \in \ZZ$,
	$$\sum_{ \substack{s,t\in\ZZ,\ s+t=k}}\widehat{g}_{s,t}=0,$$
	and thus, $\sum_{q\in\mathbb{Z}}\widehat{g}_{m-q,q-n}=0$. It follows that
	\begin{align}\notag
		\notag|I_1(m,n)|&=\bigg|\sum_{q\in\mathbb{Z}\backslash\{0\}}\widehat{g}_{m-q,q-n}S_{q,q}\bigg|
		=\bigg|\sum_{q\in\mathbb{Z}\backslash\{0\}}\frac{1}{2 q_*}\,\widehat{g}_{m-q,q-n}
		-\frac{1}{2m_*}\sum_{q\in\mathbb{Z}}\,\widehat{g}_{m-q,q-n}
		\bigg|\\
		\notag \label{eq:part2to1:subtrzero}&\le \frac{1}{2}\sum_{q\in\mathbb{Z}}\frac{(m-q)_*}{q_* m_*}\, |\widehat{g}_{m-q,q-n}|,\quad\mbox{for all } m,n\in\mathbb{Z},
	\end{align}
	and therefore,
	\begin{align*}
		J_1^{(1)}
		&\leq C\sum_{|m|>N}\bigg(\sum_{n\in\ZZ}|\beta_n|\sum_{q\in\mathbb{Z}}\frac{(m-q)_*}{q_* m_*}\,\frac{1}{(m-q)_*^{3}}\,\frac{1}{(q-n)_*^{2}}\bigg)^2\\
		&\leq  C\sum_{|m|>N}\frac{1}{m^2}\bigg(\sum_{n\in\ZZ}|\beta_n|\sum_{q\in\mathbb{Z}}\frac{1}{q_* }\,\frac{1}{(m-q)_*^{2}}\,\frac{1}{(q-n)_*^{2}}\bigg)^2.
	\end{align*}
	Since $\sup_{n}|\beta_n|\le 1$,
	$\sum_{n\in\mathbb{Z}}\frac{|\beta_n|}{(q-n)_*^{2}}\le \sum_{n\in\mathbb{Z}}\frac{1}{(q-n)_*^{2}}=O(1).
	$
	This leads to 
	\begin{align}\notag
		J_1^{(1)}&\leq C\sum_{|m|>N}\frac{1}{m^2}\bigg(\sum_{q\in\ZZ}\frac{1}{q_*(m-q)_*^{2}}\bigg)^2
		\leq C\sum_{|m|>N}\frac{1}{m^2}\left(\frac{1}{|m|^{\min\{1,2\}}}\right)^2 \leq \frac{C}{N^3},
	\end{align}
	where Lemma \ref{firstlemma} is used for the second inequality.
	
	We next consider $J_1^{(2)}$.
	Applying the Cauchy--Schwarz inequality, we obtain
	\begin{align} \notag
		\sum_{|n|>N}\frac{|\beta_n|}{(q-n)_*^{2}}\le\Big(\sum_{|n|>N}(q-n)_*^{-4}\Big)^{1/2}\leq  C f(q,N), \quad
		f(q,N)=
		\begin{cases}
			1&\mbox{if }|q|>N,\\
			(N+1-|q|)^{-3/2}&\mbox{if }|q|\le N.
		\end{cases}
	\end{align}
	As before,
	\begin{align*}
		J_1^{(2)}
		&\leq C  \sum_{|m|\leq N}\frac{1}{m_*^2}\Big(\sum_{|q|\leq N}\frac{f(q,N)}{q_* (m-q)_*^{2}}\Big)^2 +   C\sum_{|m|\leq N}\frac{1}{m_*^2}\Big(\sum_{|q|> N}\frac{f(q,N)}{q_* (m-q)_*^{2}}\Big)^2 =: C(A+B).
	\end{align*}
	Another application of the Cauchy--Schwarz inequality and Lemma \ref{firstlemma} give
	\begin{align*}
		A&=\sum_{|m|\leq \frac{N}{2}} \frac{1}{m_*^2}\Big(\sum_{|q|\leq N}\frac{f(q,N)}{q_* (m-q)_*^{2}}\Big)^2
		+\sum_{\frac{N}{2}<|m|\leq N} \frac{1}{m_*^2}\Big(\sum_{|q|\leq N}\frac{f(q,N)}{q_* (m-q)_*^{2}}\Big)^2\\
		&\le \sum_{|m|\le \frac{N}{2}} \frac{1}{m_*^2}\sum_{|q|\le N}\frac{1}{q_*^2}\Big(\sum_{|q|\le N}\frac{f(q,N)^2}{(m-q)_*^{4}}\Big)
		+\sum_{\frac{N}{2}<|m|\le N} \frac{1}{m_*^2}\Big(\sum_{|q|\le N}\frac{1}{q_*^2 (m-q)_*^{2}}\Big)\Big(\sum_{|q|\le N}\frac{f(q,N)^2}{(m-q)_*^2}\Big),
	\end{align*}
	and
	\begin{align*}
		A&\leq C\sum_{|m|\le \frac{N}{2}}\frac{1}{m_*^2}\frac{1}{(N+1-|m|)^{3}}
		+ C\sum_{\frac{N}{2}<|m|\le N}\frac{1}{m_*^2}\frac{1}{m_*^2}=O \left(\frac{1}{N^{3}} \right)+O \left(\frac{1}{N^{3}}\right).
	\end{align*}
	Similarly, 
	\begin{align*}
		B=\sum_{|m|\leq \frac{N}{2}}\frac{1}{m_*^2} \bigg(\sum_{|q|> N}\frac{f(q,N)}{q_* (m-q)_*^{2}}\bigg)^2
		+\sum_{\frac{N}{2}<|m|\leq N}\frac{1}{m_*^2} \bigg(\sum_{|q|> N}\frac{f(q,N)}{q_* (m-q)_*^{2}}\bigg)^2=O \left(\frac{1}{N^{4}}\right)+O \left(\frac{1}{N^{3}}\right),
	\end{align*}
	and therefore, $J_1^{(2)}=O \left(\frac{1}{N^{3}}\right).$
	
	\smallskip
	\smallskip

	\noindent{\underline{\bf{Estimates for $J_2^{(1)}$} and $J_2^{(2)}$.}} 
	From \eqnref{S} we see that
	\begin{align*}
		I_2(m,n)&=\sum_{s,t\in\ZZ,\ t+n=m-s=0} \widehat{g}_{s,t}\, S_{t+n,m-s}
		=\widehat{g}_{m,-n}\, \log \gamma.
	\end{align*}
	The Cauchy--Schwarz inequality is again applied and one arrives at the estimates
	\begin{align*}
		J_2^{(1)}
		&=(\log\gamma)^2\sum_{|m|>N}\bigg|\sum_{n\in\ZZ} \beta_n\, \widehat{g}_{m,-n} \bigg|^2
		\leq (\log\gamma)^2\sum_{|m|>N} \sum_{n\in\ZZ}\left|\widehat{g}_{m,-n}\right|^2
		\leq C \sum_{|m|>N}\sum_{n\in\ZZ}\frac{1}{m_*^{6} n_*^{4}},\\
		J_2^{(2)}
		&=(\log \gamma)^2\sum_{|m|\leq N}\bigg|\sum_{|n|>N} \beta_n\, \widehat{g}_{m,-n}\bigg|^2
		\leq (\log \gamma)^2 \sum_{|m|\leq N}\sum_{|n|>N}\left|\widehat{g}_{m,-n}\right|^2
		\leq C \sum_{|m|\leq N}\sum_{|n|>N}\frac{1}{m_*^{6} n_*^{4}},
	\end{align*}
	and hence, $J_2^{(1)}, J_2^{(2)}=O\left(\frac{1}{N^{3}}\right)$. 
	
	\smallskip
	\smallskip
	
	\noindent{\underline{\bf{Estimates for $J_3^{(1)}$} and $J_3^{(2)}$}.} 
	The inequality  \eqnref{S2} will be the main tool of the discussion here.
	Fix $m\in\ZZ$. Then
	\begin{align}\notag
		&\Big|\sum_{n\in\ZZ} \beta_n I_3(m,n)\Big|^2
		=\Big|\sum_{n\in\ZZ}\sum_{\substack{s,t\in\ZZ \\ (t+n)(m-s)<0}} \beta_n \,\widehat{g}_{s,t}\,S_{n+t,m-s}\Big|^2\leq CAB,
	\end{align}
	for some constant $C$ and by Lemma \ref{firstlemma} and the fact that $|\beta_n|\leq 1$, we have
	\begin{align*}
		A&= 
		\sum_{n\in\ZZ}\sum_{\substack{s,t\in\ZZ \\ (t+n)(m-s)<0}} |\beta_n|\frac{(n+t)_*^{-2}}{s_*^{2} \,t_*^{2}}
		\leq\sum_{n\in\ZZ}\sum_{t\in\ZZ}\frac{(n+t)_*^{-2}}{t_*^{2}}
		\sum_{s\in\ZZ}\frac{1}{s_*^{2}}
		\leq C\sum_{n\in\ZZ}\frac{1}{n_*^2}=O(1),\\
		B&= 
		\sum_{n\in\ZZ}\sum_{\substack{s,t\in\ZZ\\ (t+n)(m-s)<0}} |\beta_n|\,\frac{(n+t)_*^2\,|S_{n+t,m-s}|^2}{s_*^{4}\, t_*^{2}}
		\leq C\sum_{n\in\ZZ}\sum_{s,t\in\ZZ}\frac{(n+t)_*^2\, (n+t)_*^{-2-2\alpha}(m-s)_*^{-4}}{s_*^{4}\, t_*^{2}}\\
		&\leq C\sum_{n\in\ZZ}\sum_{t\in\ZZ}\frac{(n+t)_*^{-2\alpha}}{t_*^2}\sum_{s\in\ZZ}\frac{(m-s)_*^{-4}}{s_*^{4}}
		\leq C\sum_{n\in\ZZ}\frac{1}{n_*^{2\alpha}}\frac{1}{m_*^{4}}
		= O\left(\frac{1}{m_*^{4}}\right).
	\end{align*}
	Therefore,
	$
	J_3^{(1)}\leq C \sum_{|m|>N} AB\leq C \sum_{|m|>N} \frac{1}{m_*^{4}}=O \left(\frac{1}{N^{3}}\right).
	$

	\smallskip
	
	Proceeding in a similar fashion as before,
	\begin{align}\notag
		&\Big|\sum_{|n|>N}\beta_n I_3(m,n)\Big|^2
		=\Big|\sum_{|n|>N}\sum_{\substack{s,t\in\ZZ \\ (t+n)(m-s)<0}} \beta_n \,\widehat{g}_{s,t}\,S_{n+t,m-s}\Big|^2\leq C\widetilde{A}\widetilde{B},
	\end{align}
	where
	\begin{align*}
		\widetilde{A}&=  \sum_{|n|>N}\sum_{\substack{s,t\in\ZZ \\ (t+n)(m-s)<0}} |\beta_n|\frac{(n+t)_*^{-2}}{s_*^{3} \,t_*^{2}}
		\leq\sum_{|n|>N}|\beta_n|\sum_{t\in\ZZ}\frac{(n+t)_*^{-2}}{t_*^{2}}
		\sum_{s\in\ZZ}\frac{1}{s_*^{3}}\\
		&\leq C\sum_{|n|>N}|\beta_n|\,\frac{1}{n_*^2}\leq C\Big(\sum_{|n|>N}|\beta_n|^2\Big)^{\frac{1}{2}}\Big(\sum_{|n|>N}\frac{1}{n_*^4}\Big)^{\frac{1}{2}} =O(N^{-\frac{3}{2}}),\\
		\widetilde{B}&= \sum_{|n|>N}\sum_{\substack{s,t\in\ZZ\\ (t+n)(m-s)<0}} |\beta_n|\,\frac{(n+t)_*^2\,|S_{n+t,m-s}|^2}{s_*^{3}\,t_*^{2}}\\
		&\leq C\sum_{|n|>N}|\beta_n|\sum_{t\in\ZZ}\frac{(n+t)_*^{-1-2\alpha}}{t_*^2}\,\sum_{s\in\ZZ}\frac{(m-s)_*^{-3}}{s_*^{3}}
		\leq C\sum_{|n|>N}|\beta_n|\frac{1}{n_*^{2}}\frac{1}{m_*^{3}}
		= O(m_*^{-3}N^{-\frac{3}{2}}),
	\end{align*}

	%
	and therefore,
	$J_3^{(2)}\leq C\sum_{|m|\leq N}\widetilde{A}\widetilde{B}
	\leq C \sum_{|m|\leq N} m_*^{-3} N^{-3}=O(N^{-3}).$

	\smallskip
	
	Since $J_{k}^{(1)}, J_k^{(2)} \leq \frac{C}{N^2}$ for all $k=1,2,3,$ and for some constant $C$ independent of $\beta_n$ and $N$, the result follows from inequality (\ref{Js}).
	
\end{proof}

\subsubsection{Proof of Theorem \ref{theorem:main:smooth}}

\begin{proof}
	Let $\varphi=\sum_{n\in\ZZ}\beta_n \widetilde{\eta}_n \in L^2(\p\Om)$ with $\|\varphi\|^2_0=\sum_{n\in\ZZ}|\beta_n|^2=1$.
	Recall first from (\ref{Kcaldecomp}) that 
	$$
	\Kcal_{\Om}^\omega=\Kcal_{\Om}^0 + \mathcal{R}_1+\mathcal{R}_2, 
	$$
	where by Theorem \ref{lem:uss}, 
	\begin{align*}
		\mathcal{K}_\Om^0[\varphi]&=\frac{1}{2}\beta_0\widetilde{\eta}_0+\frac{1}{2}\sum_{k=1}^\infty \left( \sum_{n=1}^\infty \beta_{n} \frac{c_{n,k}}{\gamma^{n+k}}\right)\widetilde{\eta}_{-k}+\frac{1}{2}\sum_{k=1}^\infty \left(\sum_{n=1}^\infty \beta_{-n} \frac{\overline{c_{n,k}}}{{\gamma^{n+k}}}\right)\widetilde{\eta}_{k},
	\end{align*}
	using the Grunsky coefficients $c_{n,k}$. Then the Cauchy--Schwarz inequality, the assumption that $\alpha>\frac{1}{2}$, \eqnref{Grun:inter}, and Theorem \ref{Suetin} together give 
	\begin{align}\notag
		\|(\mathcal{K}_\Om^0-P_N\mathcal{K}_\Om^0 P_N)[\varphi]\|_{0}^2
		=&\frac{1}{4}\sum_{k=1}^\infty \bigg|\sum_{n:\max\{k,n\}>N}\beta_n\frac{c_{n,k}}{\gamma^{n+k}}\bigg|^2 +\frac{1}{4}\sum_{k=1}^\infty \bigg|\sum_{n:\max\{k,n\}>N}\beta_{-n}  \frac{\overline{c_{n,k}}}{{\gamma^{n+k}}}\bigg|^2\\\notag
		\le&\frac{1}{4}\sum_{k=1}^\infty\, \sum_{n:\max\{k,n\}>N}\left|\frac{c_{n,k}}{\gamma^{n+k}}\right|^2 +\frac{1}{4}\sum_{k=1}^\infty\, \sum_{n:\max\{k,n\}>N}\left|  \frac{\overline{c_{n,k}}}{{\gamma^{n+k}}}\right|^2
		\\\notag
		=&\frac{1}{2}\sum_{k=1}^\infty\, \sum_{n:\max\{k,n\}>N}\frac{k}{n}\frac{c_{n,k}\overline{c_{k,n}}}{\gamma^{2n+2k}}
		\leq\frac{1}{2N^3}\sum_{k,n=1}^\infty\frac{k}{n}\frac{M}{\max\{k,n\}^{1+2\alpha}}=O\left(\frac{1}{N^3}\right).
	\end{align}
	It follows that
	\beq\label{eq1}
	\left\|\mathcal{K}_\Om^0-P_N \mathcal{K}_\Om^0 P_N\right\|_{0}^2\leq M N^{-3},
	\eeq
	for some constant $M$ independent of $N$ and $\omega$. 
	
	For the term involving $\mathcal{R}_1$, the Cauchy--Schwarz inequality and Lemma \ref{R1} imply that
	\beq\label{eq2}
	\begin{aligned}
		\left\|\mathcal{R}_1-P_N \mathcal{R}_1 P_N\right\|_{0}^2\leq & \sum_{\max\{|m|,|n|\}>N}|\langle \mathcal{R}_1[\widetilde{\eta}_n], \widetilde{\eta}_m\rangle_0|^2 \leq \sum_{\max\{|m|,|n|\}>N} C(w)^2(1+|m|)^{-6}(1+|n|)^{-4} \\
		& \leq \widetilde{C}(\omega)N^{-3},
	\end{aligned}
	\eeq
	for $m,n \in \mathbb{Z}$ and $C(\omega)$ and $\widetilde{C}(\omega)$ that depend continuously on $\omega \in \mathbb{C} \backslash \sqrt{-1} \mathbb{R}^{-}$.
	Lastly, we know from Lemma \ref{R2estimate} that
	\beq\label{eq3}
	\|\mathcal{R}_2-P_N\mathcal{R}_2 P_N\|_{0}^2=O({N^{-3}}).
	\eeq
	Combining \eqnref{Kcaldecomp}, \eqnref{eq1}, \eqnref{eq2}, and \eqnref{eq3} proves the result.
\end{proof}

\subsection{Matrix representation of $P_n(T_n\mathcal{K}^{\omega}_{\Omega})P_n$
}
\label{subsection:K:matrix}

For each $m \in \mathbb{Z}$, we expand the operators given by \eqnref{Kcal_k_exp:2} in $\{\teta_k\}$ as
\beq\label{def:dkm}
\begin{aligned} \notag
	\Kcal^{(j,s)}_{\Om}\left[\teta_{ m}\right]
	=\sum\nolimits_{k\in\NN} d^{(j,s)}_{k,m}\,\teta_k \quad \mbox{for } j\geq 0,\ s=1,2,
\end{aligned}
\eeq
for some $d^{(j,s)}_{k,m} \in \mathbb{C}$. Then
\begin{align}\notag
	P_n\left(T_n \Kcal_\Om^\omega\right)P_n\left[\teta_{ m}\right]
	=&
	\begin{cases}
		\ds \sum\nolimits_{|k|\leq n} d_{k,m} \teta_k \quad &\mbox{if }|m|\leq n,\\
		\ds 0,\quad &\mbox{if }|m|>n,
	\end{cases}
\end{align}
with
\beq\label{dkm_def} \notag
d_{k,m}=\sum\nolimits_{0\leq j\leq n} \left(\omega^{2j} \log \omega\right)d_{k,m}^{(j,1)}+ \sum\nolimits_{0\leq j\leq n} \omega^{2j} \,d_{k,m}^{(j,2)}.
\eeq
Hence, one can identify $P_n\left(T_n \Kcal_\Om^\omega\right)P_n$ with the finite-dimensional complex matrix
\beq\label{mat:finite:K} \notag
\mathbb{K}_{\Om,n}^\omega := \left(d_{k,m}\right)_{-n \leq k,m \leq n}.
\eeq

\begin{theorem} \label{lemma:dlp}
	Let $\Psi$ be a finite Laurent series. For a fixed $n$, $d_{k,m}^{(j,1)}$ and $d_{k,m}^{(j,2)}$ with $0\leq j\leq n$ can be explicitly computed using a finite number of operations of scalar additions and multiplications on $\gamma$, $\{a_n\}_{n=0}^\infty$, their reciprocals, conjugates, and other constants.
\end{theorem}

\begin{proof}
	For $z, \tau \in \mathbb{C}$ and $j \in \mathbb{N}$, we have the binomial expansion
	\begin{equation}\label{eq:ln:taylor2} \notag \begin{aligned}
			|z-\tau|^{2j}&=\sum_{j_1=0}^j\sum_{j_2=0}^j{j\choose j_1}{j\choose j_2} (-1)^{j_1+j_2}\, z^{j-j_1}\,\tau^{j_1}\,\overline{z^{j-j_2} \,\tau^{j_2}}.\end{aligned}
	\end{equation}
	Using \eqnref{eq:fundamental:sol} and the fact that $\Gamma^0(z)=(1/2\pi)\log|z|$, we see that 
	\begin{equation}\label{eq:3p3:kernel}
		\begin{aligned}
			&(\Gamma^\omega-\Gamma^0)(z-\tau)\\
			=&\tau(\omega)+\sum_{j=1}^\infty  \left(b_j \log (\omega|z-\tau|) +c_j \right) \omega^{2j}
			\sum_{j_1=0}^j\sum_{j_2=0}^j {j\choose j_1}{j\choose j_2}(-1)^{j_1+j_2}\,z^{j-j_1}\,\tau^{j_1}\,\overline{z^{j-j_2}}\,\overline{\tau^{j_2}}.
		\end{aligned}
	\end{equation}
	
	The integral kernel of $\mathcal{K}_\Omega^\omega-\mathcal{K}_\Omega^0$ is the term-by-term normal derivative of \eqnref{eq:3p3:kernel} in $\tau$. Then for $\varphi\in L^2(\partial\Omega)$ and $z\in\partial{\Omega}$,
	\begin{align}\notag
		\Kcal_\Om^{(j,1)}[\varphi](z)&=b_j \sum_{j_1=0}^j\sum_{j_2=0}^j{j\choose j_1}{j\choose j_2} (-1)^{j_1+j_2} z^{j-j_1}\overline{z^{j-j_2}}\int_{\p\Omega}\frac{\p (\tau^{j_1}\,\overline{\tau^{j_2}})}{\p\nu_\tau}\,\varphi(\tau)\,d\sigma(\tau),\\ \notag
		\Kcal_\Om^{(j,2)}[\varphi](z)&=\sum_{j_1=0}^j\sum_{j_2=0}^j{j\choose j_1}{j\choose j_2} (-1)^{j_1+j_2} z^{j-j_1}\overline{z^{j-j_2}}
		\Bigg[ 2\pi b_j \,\mathcal{K}_\Omega^0\big[\tau^{j_1}\,\overline{\tau^{j_2}}\,\varphi(\tau)\big](z)\\ \notag
		&+2\pi b_j \,\mathcal{S}_\Omega^0\bigg[\frac{\p  (\tau^{j_1}\,\overline{\tau^{j_2}})}{\p\nu_\tau}\,\varphi(\tau)\bigg](z)
		+ c_j\int_{\p\Omega}\frac{\p (\tau^{j_1}\,\overline{\tau^{j_2}})}{\p\nu_\tau}\,\varphi(\tau)\,d\sigma(\tau)\Bigg].
	\end{align}
	By \eqnref{pOm:normalderi} with $\tau=\Psi(\xi)$ and  $\xi=e^{\rho+\mathrm{i} \theta}$,
	\begin{align} \notag
		h(\tau)\,\frac{\p (\tau^{j_1}\overline{\tau^{j_2}})}{\p\nu_\tau}\bigg|_{\p\Om}^{+}(\tau)
		=&\lim_{\rho\to(\log \gamma)+} \frac{\p}{\p \rho} \left(\tau^{j_1}\overline{\tau^{j_2}}\right)\\[1mm]
		\notag
		=&j_1 \Psi(\xi)^{j_1-1}\,\Psi'(\xi)\,\xi\,\overline{\Psi(\xi)^{j_2}}
		+j_2\Psi(\xi)^{j_1}\,\overline{\Psi(\xi)^{j_2-1}\,\Psi'(\xi)\,\xi}, 
	\end{align}
	and using $d\sigma(\tau)=h(\tau)d\theta$, we obtain
	\begin{align}\notag
		&\int_{\p\Omega}\frac{\p (\tau^{j_1}\,\overline{\tau^{j_2}})}{\p\nu_\tau}\, \varphi (\tau)\,d\sigma(\tau)\\ 
		=&\int_0^{2\pi} \left(j_1 \Psi(\xi)^{j_1-1}\,\Psi'(\xi)\,\xi\,\overline{\Psi(\xi)^{j_2}}
		+j_2\Psi(\xi)^{j_1}\,\overline{\Psi(\xi)^{j_2-1}\,\Psi'(\xi)\,\xi}\right)\varphi(\Psi(\gamma e^{\mathrm{i}\theta}))\,d\theta.\label{comp:part_Kcal}
	\end{align}
	
	If $\Psi(\xi)$ is a finite Laurent series, then setting $\varphi=\eta_m$ (for a fixed $m$), \eqnref{comp:part_Kcal} can be calculated using a finite number of additions and multiplications of complex numbers. This also holds for the expansion of $\tau^{j_1}\overline{\tau^{j_2}}\eta_m(\tau)$ in terms of the density functions $\{\eta_k(\tau)\}_{k\in\ZZ}$ and for the expansion of $\frac{\p  (\tau^{j_1}\,\overline{\tau^{j_2}})}{\p \nu_\tau}\big|_{\p\Om}^{+}\,\eta_m(\tau)$ in terms of $\{\zeta_k(\tau)\}_{k\in\ZZ}$. Furthermore, for fixed $j_1$ and $j_2$, $z^{j-j_1}\overline{z^{j-j_2}}$ admits an expansion into $\{\eta_k(z)\}_{k\in\ZZ}$ by finitely many elementary operations of complex numbers. 
\end{proof}

\section{Approximation of Laplacian Eigenvalues}\label{subsec:multiple}

We can estimate the Laplacian eigenvalues by investigating the operator $\pm \frac{1}{2}+\mathcal{K}_\Om^\omega$ using the following result of M. Mitrea:
\begin{theorem}[{\cite[Theorems 8.6, 8.7, 8.9]{Mitrea:1996:BHH}}]\label{Lemma:eigen:DN}
	\label{lem:characterize:eig} Let $\Omega$ be a bounded simply connected $C^{1,\alpha}$ domain in $\mathbb{R}^2$ for some $\alpha>0$. For $\omega_0>0$,
	\begin{itemize}
		\item[\rm(a)] $\omega_0^2$ is a Neumann eigenvalue of $-\Delta$ on $\Om$ if and only if $-\frac{1}{2}I+\mathcal{K}_\Om^{\omega_0}$ is not invertible on $L^2(\p\Om)$.
		\item[\rm(b)] $\omega^2_0$ is a Dirichlet eigenvalue of $-\Delta$ on $\Om$ if and only if $\frac{1}{2}I+\mathcal{K}_\Om^{\omega_0}$ is not invertible on $L^2(\p\Om)$.
	\end{itemize}
\end{theorem}

In this section, we apply the Gohberg--Sigal Theory to approximate the Neumann eigenvalues of the Laplacian using the characteristic values of the operator $-\frac{1}{2}I+\mathcal{K}_{\Omega}^{\omega}$. 
All the results in the section can be proved analogously for the Dirichlet eigenvalues by considering $\frac{1}{2}I+\mathcal{K}_{\Omega}^{\omega}$ instead. 

For each Neumann eigenvalue $\omega_0^2 \neq 0$ $(\omega_0 > 0)$ of the Laplacian, there exists a $\delta_{\omega_0}>0$ such that for every $0<\delta<\delta_{\omega_0}
$, the disk $B_{\delta}(\omega_0)$ centered at $\omega_0$ is compactly contained in $\mathbb{C}\backslash\sqrt{-1}\mathbb{R}^-$ and
$$-\frac{1}{2}I+\Kcal_\Om^\omega\mbox{ is invertible on }L^2(\p\Om)\mbox{ for all }w \in B_{\delta}(\omega_0)  \backslash \{\omega_0\}.$$

\subsection{Application of the Gohberg--Sigal Theory to finite section operators of $\Kcal_\Om^\omega$}

The first lemma can be proved similarly as in Lemma 6.1 and Theorem 6.1 of \cite{Ammari:2004:SRS} by H. Ammari and F. Triki. See also \cite[Theorem 2] {Steinberg:1968:MFC}. 
\begin{lemma} \label{AT}
	A Neumann eigenvalue $\omega_0^2 \neq 0$ of $-\Delta$ is a simple pole of $\Big(-\frac{1}{2}I+\Kcal_\Om^\omega\Big)^{-1}$ and in a small neighborhood of $\omega_0$, one can write 
	$$
	\Big(-\frac{1}{2}I+\Kcal_\Om^\omega\Big)^{-1}=\frac{1}{
		\omega-\omega_0}\mathcal{Q}+\mathcal{R}(\omega),
	$$
	for some rank $\mu$ operator $\mathcal{Q}:  \ker \Big(-\frac{1}{2}I+\Kcal_\Om^{\omega_0}\Big)^* \rightarrow \ker \Big(-\frac{1}{2}I+\Kcal_\Om^{\omega_0}\Big)$
	and an analytic family $\mathcal{R}(\omega)$ of bounded linear operators.
\end{lemma}
\begin{lemma}
	\label{lemma:upperboundofTrace} \notag
	Let $\p\Om$ be of class $C^{1,\alpha}$ for some $0<\alpha<1$. Let $\omega_0^2 \neq 0$ be a Neumann eigenvalue of $-\Delta$ of multiplicity $\mu$. Fix $0<\delta<\delta_{\omega_0}$ and consider an analytic family of bounded linear operators $X(\omega)$ on $L^2(\p\Om)$ for $\omega\in B_{\delta}(\omega_0)$. Then for each $p\ge m \geq 1$, 
	\begin{align}
		\notag&\left|\operatorname{tr}\int_{\p B_{\delta}(\omega_0)}\left(\omega-\omega_0\right)^{m-1}\left[\Big(-\frac{1}{2}I+\Kcal_\Om^\omega\Big)^{-1}X(\omega)\right]^p\,d\omega\right|\\
		\label{eq:multipleeig:trace}\le&\, 2\pi  p^2 \mu \,\delta^m \sup_{\omega\in\p B_{\delta}(\omega_0)}\left\|\Big(-\frac{1}{2}I+\Kcal_\Om^\omega\Big)^{-1}X(\omega)\right\|_{\mathcal{L}(L^2(\p\Om))}^p.
	\end{align}
	In addition, the left hand side of the inequality is zero for all $p=0,1,\dots,m-1$.
\end{lemma}
\begin{proof}
	
	We note first that $\omega_0$ is a characteristic value of 
	$
	A(\omega)=-\frac{1}{2}I+\Kcal_\Om^\omega
	$
	with
	\beq\label{eq:simplepole:multipleeig:0} \notag
	M(A(\omega_0))=\dim\ker A(\omega_0) = \mu.
	\eeq
	It follows from Lemma \ref{AT} that
	$$
	\label{eq:simplepole:multipleeig}
	A(\omega)^{-1}=\frac{1}{
		\omega-\omega_0}\mathcal{Q}+\mathcal{R}(\omega),
	$$
	for some rank $\mu$ operator $\mathcal{Q}: \text{ ran }A(\omega_0)^{\perp} \rightarrow \ker A(\omega_0)$
	and an analytic family $\mathcal{R}(\omega)$ of bounded linear operators.

	We next expand the analytic families $X(\omega)$ and $\left[\mathcal{R}(\omega)X(\omega)\right]^k$ for $k \in \mathbb{N}$ as Laurent series about $\omega=\omega_0$ so that 
	$$
	X(\omega)=\sum_{j\ge 0} \left(\omega-\omega_0\right)^{j} X_j,\quad\left[\mathcal{R}(\omega)X(\omega)\right]^k=\sum_{j\ge 0} \left(\omega-\omega_0\right)^{j} R_{kj},
	$$
	for some operators $X_j$, $R_{kj}$ on $L^2(\p\Om)$.
	It follows that 
	\begin{align*}
		\left[A(\omega)^{-1} X(\omega)\right]^{p}&=\left[(\omega-\omega_0)^{-1}\mathcal{L}X(\omega)+\mathcal{R}(\omega) X(\omega)\right]^{p}\\
		&=\left[\left((\omega-\omega_0)^{-1}\mathcal{L}X_0 + \mathcal{L}X_1 + (\omega-\omega_0)\mathcal{L}X_0 +\cdots\right) + \mathcal{R}(\omega) X(\omega)\right]^{p}.
	\end{align*}
	We note that for $m=1,\dots,p$, the operator coefficient of the term $(\omega-\omega_0)^{-m}$ equals $$\mathcal{L}Y_{00m}+\sum_{k=1}^{p-1}\sum_{j=0}^{p-m} {R}_{kj}\mathcal{L}Y_{kjm},$$ where the operators $Y_{kjm}$ are independent of $\omega$.
	Moreover, when $m >p$, this coefficient is zero. Then by the Residue Theorem for operator-valued functions, 
	\begin{align}\label{eq:resthm:multipleeig}
		\int_{\p B_\delta(\omega_0)}\left(\omega-\omega_0\right)^{m-1}\left[A(\omega)^{-1}X(\omega)\right]^p\,d\omega =
		\begin{cases}
			\ds 2\pi{\rm i}\bigg(\mathcal{L}Y_{00m}+\sum_{k=1}^{p-1}\sum_{j=0}^{p-m} {R}_{kj}\mathcal{L}Y_{kjm}\bigg) &\mbox{if }p\ge m,\\
			\ds 0&\mbox{if }p<m.
		\end{cases}
	\end{align}
	This proves the lemma for the case $0 \leq p<m$.
	
	For $p\ge m$, since $\text{ran }\mathcal{L} = \ker A(\omega_0)$ is of dimension $\mu$, denote by $\{v_1,\dots,v_\mu\} \subset L^2(\p\Om)$ the set of basis elements of $\mathcal{L}$. The right-hand side of \eqnref{eq:resthm:multipleeig} is contained in the finite dimensional space
	$$S:=\operatorname{span}\Big(\{v_1,\dots,v_\mu\}\cup\{R_{kj}v_l\,:\,0\le j\le p-m,\ 1\le k\le p-1,\ 1\le l\le \mu\}\Big),$$
	with $\dim S \le p^2\mu$.
	Denote by $\{v_k^*\}_{k=1}^{\dim S}$ an orthonormal basis for $S$.
	Using the trace of a linear operator, we then have
	\begin{align*}
		&\operatorname{tr}\int_{\p B_\delta(\omega_0)}\left(\omega-\omega_0\right)^{m-1}\left[A(\omega)^{-1}X(\omega)\right]^p\,d\omega\\
		=&\sum_{k=1,\dots,\dim S}\Big\langle v_k^*,\int_{\p B_{\delta}(\omega)}\left(\omega-\omega_0\right)^{m-1}\left[A(\omega)^{-1}X(\omega)\right]^p v_k^*\,d\omega \Big\rangle_{L^2(\p\Om)}\\
		=&\sum_{k=1,\dots,\dim S}\int_{\p B_\delta(\omega_0)}\left(\omega-\omega_0\right)^{m-1}\langle v_k^*,\left[A(\omega)^{-1}X(\omega)\right]^p v_k^*\rangle_{L^2(\p\Om)}\,d\omega.
	\end{align*}
	We complete the proof using the Cauchy-Schwarz inequality.
\end{proof}
\begin{lemma}\label{lemma:inequality:sympwrsum} Fix $\mu\in\NN$ and $R>0$. Let $\{x_m\}_{m=1}^\mu \subset \mathbb{C}$ satisfy
	$\left|\sum_{m=1}^\mu x_m^n\right| \le R^n$ for all $n=1,\cdots,\mu.$
	Then $|x_m|\le 2R$ for all $m=1,\cdots,\mu$.
\end{lemma}
\begin{proof}
	For each $n=1,2,\cdots,\mu$, we set
	$$y_n:=\sum_{m=1}^\mu x_m^n\quad\mbox{and}\quad c_n:=(-1)^n\sum\bigg\{\prod_{m=1}^\mu x_m^{a_m}\,:\,a_m\in\{0,1\}\mbox{ and }\sum_{m=1}^\mu a_m=n\bigg\}.$$
	Since
	$
	\prod_{m=1}^\mu(x-x_m)=x^\mu+\sum_{n=1}^\mu c_n x^{\mu-n},
	$
	the inequality (see, for example, \cite[Chapter VI.2]{Yap:2000:FPA}) 
	\beq\label{eq:LgrgKnth:ineq} \notag
	\sup_{1\le m\le\mu}|x_m|\le 2\max\left\{\left|c_n\right|^{1/n}\,:\,1\le n\le \mu\right\}\eeq
	holds. It is also easy to check the recursion relation
	\beq\notag
	c_1=-y_1\quad\mbox{and}\quad c_n=-\frac{1}{n}\bigg(y_n+\sum_{m=1}^{n-1} y_{n-m}c_m\bigg),\quad\mbox{for }n=2,3,\cdots,\mu,
	\eeq
	using a simple combinatorial argument. One can then use induction on $n$ to show that the assumption $|y_n|\le R^n$ implies that $|c_n|\le R^n$ for $n=1,2,\cdots,\mu$.
\end{proof}

For the following results, let
$A(\omega)=-\frac{1}{2}I+\Kcal_\Om^\omega$ and $A_n(\omega)=-\frac{1}{2}I+T_n\Kcal_\Om^\omega$. The operator $A(\omega)^{-1}$ is bounded since $\partial \Omega$ is of $C^{1, \alpha}$.

\begin{theorem} \label{ineq:rho:parta}
	We have
	\beq \notag
	\sup_{\omega\in\p B_\delta(\omega_0)}\left\|A(\omega)^{-1}\left(A(\omega)-A_n(\omega)\right)\right\|_{\mathcal{L}(L^2(\p\Om))} \le \frac{C}{n!}.
	\eeq   
\end{theorem}
\begin{proof}
	In the proof of Lemma \ref{lemma:kminustk:pointwise}, the positive numbers $C$ and $K$ in inequality \eqref{inequality:KmTK} are constants independent of $\omega\in\mathbb{C}\backslash\sqrt{-1}\mathbb{R}^{-}$. Therefore, for some constant $C>0$, we have
	\beq\label{ineq:smallcircle:factorialdecay}
	\sup_{\omega\in\p B_\delta(\omega_0)}\left\|A(\omega)-A_n(\omega)\right\|_{\mathcal{L}(L^2(\p\Om))} \le \frac{C}{n!}.
	\eeq
	Using Lemma \ref{AT}, we obtain
	\beq\label{ineq:smallcircle:finite}
	\sup_{\omega\in\p B_\delta(\omega_0)}\left\|A(\omega)^{-1}\right\|_{\mathcal{L}(L^2(\p\Om))}\le \frac{1}{\delta}\left\|\mathcal{Q}\right\|_{\mathcal{L}(L^2(\p\Om))} + \sup_{\omega\in\p B_\delta(\omega_0)}\left\|\mathcal{R}(\omega)\right\|_{\mathcal{L}(L^2(\p\Om))}<\infty,
	\eeq
	where $\mathcal{Q}$ is extended to be a bounded operator on $L^2(\p\Om)$ by assigning $0$ on the orthogonal complement of $\ker A(\omega)^*$. Here, the upper bound is finite because $\mathcal{Q}$ is bounded and $\mathcal{R}(\omega)$ is analytic.
	Combining \eqref{ineq:smallcircle:factorialdecay} and \eqref{ineq:smallcircle:finite} gives the desired inequality.
\end{proof}

\subsection{Proofs of Theorems \ref{theorem:main:first} and \ref{theorem:main:second}}
\begin{proof}[Proof of Theorem \ref{theorem:main:first}]
	Choose $N \in \mathbb{N}$ satisfying $\frac{C}{N!}<1$. Then the operator generalization of Rouche's result given in Theorem \ref{general:rouche} shows that
	for each $n\ge N$, $A_n(\omega)= \frac{1}{2}I+T_n\Kcal_\Om^\omega$ has, $\mu$ characteristic values (counting multiplicity) in $\overline{B_\delta(\omega_0)}$.
	This proves part (a).
	
	\smallskip
	
	For part (b), we start by applying the generalized Argument Principle of Theorem \ref{arg:principle} to $f(z)=z^m$ for $m=1, \cdots, \mu,$ and $n \geq {N}$. Since 
	$$
	A_n(\omega)=A(\omega)\left[I-A(\omega)^{-1}(A(\omega)-A_n(\omega))\right],
	$$
	it follows that 
	\begin{align*} 
		\sum_{j=1}^\mu \left(\omega_{n,j}-\omega_0\right)^m 
		= &\frac{1}{2\pi{\rm i}}\operatorname{tr}\int_{\p B_\delta(\omega_0)}\left(\omega-\omega_0\right)^m A_n(\omega)^{-1}\frac{\p}{\p\omega}A_n(\omega)\,d\omega\\ 
		=&\frac{1}{2\pi{\rm i}}\sum_{p=0}^\infty \operatorname{tr}\int_{\p B_\delta(\omega_0)}\left(\omega-\omega_0\right)^m \left[A(\omega)^{-1}(A(\omega)-A_n(\omega))\right]^p A(\omega)^{-1}\frac{\p}{\p\omega}A_n(\omega)\,d\omega.
	\end{align*}
	Moreover, proceeding as in the proof of \cite[Theorem 3.9]{Ammari:2009:LPT}, one finds that
	\beq\label{eq:multifreq:argprin:2} 
	\sum_{j=1}^\mu \left(\omega_{n,j}-\omega_0\right)^m  = \frac{1}{2\pi{\rm i}}\sum_{p=1}^\infty\frac{m}{p}\operatorname{tr}\int_{\p B_\delta(\omega_0)}\left(\omega-\omega_0\right)^{m-1}\left[A(\omega)^{-1}(A(\omega)-A_n(\omega))\right]^p\,d\omega.
	\eeq
	Then by Lemma \ref{lemma:upperboundofTrace} and the fact that for $0 \leq x < 1$, $\sum_{p=m}^\infty p x^p=\frac{m-(m-1)x}{(1-x)^2}x^m$, we arrive at 
	\begin{align}\notag
		\sum_{j=1}^\mu \left(\omega_{n,j}-\omega_0\right)^m 
		&\le\sum_{p=m}^\infty \frac{m}{p}p^2\mu\delta^m{\rho_{n}}^p=m\delta^m \mu\left(\frac{m-(m-1){\rho_{n}}}{(1-\rho_{n})^2}\right){\rho_{n}}^{m} \\ \notag
		&\le \delta^m \mu^m\left(\frac{\mu^m}{(1-\rho_{n})^{2\mu}}\right){\rho_{n}}^{m} =\left(\frac{\delta\mu^2\rho_{n}}{(1-\rho_{n})^2}\right)^m,
	\end{align}
	where $\rho_n:=\sup_{\omega\in\p B_\delta(\omega_0)}\left\|A(\omega)^{-1}\left(A(\omega)-{A}_n(\omega)\right)\right\|_{\mathcal{L}(L^2(\p\Om))}$.
	
	Finally, Lemma \ref{lemma:inequality:sympwrsum} implies that for $j=1, \cdots, \mu,$
	$$ \left|\omega_{n,j}-\omega_0\right|\le \frac{2\delta\mu^2\rho_{n}}{(1-\rho_{n})^2},
	$$
	and one applies Theorem \ref{ineq:rho:parta} to complete the proof.
\end{proof}

\begin{proof}[Proof of Theorem \ref{theorem:main:second}]
	%
	We set
	$ \widetilde{A}_n(\omega)= -\frac{1}{2}I+ P_n(T_n\Kcal_{\p\Om}^\omega)P_n$
	and
	\begin{align} \notag
		\widetilde{\rho}_{n}&:=\sup_{\omega\in\p B_\delta(\omega_0)}\left\|A(\omega)^{-1}\left(A(\omega)-\widetilde{A}_n(\omega)\right)\right\|_{\mathcal{L}(L^2(\p\Om))}.
	\end{align}
	By Lemma \ref{Kconvergence} and Theorem \ref{theorem:main:smooth}, we know that 
	\beq \label{ineq:rho}
	\begin{aligned} 
		&\widetilde{\rho}_{n}\to 0\quad\mbox{as }n\rightarrow\infty,\\ 
		&\widetilde{\rho}_{n}=O\big(n^{-3/2}\big)\quad\mbox{under the assumption in (c)}.
	\end{aligned}
	\eeq
	Part (a) can be proved by choosing $N \in \mathbb{N}$ such that $\widetilde{\rho}_{n}<1$ for all $n\geq N$ and using Theorem \ref{general:rouche}.
	
	The proofs for parts (b) and (c) are similar those given in the previous theorem. For $n\geq N$ and $m=1,\dots,\mu$, one has, as in (\ref{eq:multifreq:argprin:2}), 
	\beq\label{eq:multifreq:argprin:theorem} \notag
	\sum_{j=1}^\mu \left(\tomega_{n,j}-\omega_0\right)^m  = \frac{1}{2\pi{\rm i}}\sum_{p=1}^\infty\frac{m}{p}\operatorname{tr}\int_{\p V_\delta}\left(\omega-\omega_0\right)^{m-1}\left[A(\omega)^{-1}(A(\omega)-\widetilde{A}_n(\omega))\right]^p\,d\omega,
	\eeq
	which results in
	$$\left|\tomega_{n,j}-\omega_0\right|\le \frac{2\delta\mu^2\widetilde{\rho}_{n}}{(1-\widetilde{\rho}_{n})^2},\quad\mbox{for each }j=1,\dots,\mu.$$ 
	We then recall (\ref{ineq:rho}) to complete the proof.
\end{proof}

\subsection{Approximation of Eigenfunctions}\label{subsec:eigenfunctions}

We use the notations from Theorem \ref{theorem:main:second} and consider the eigenspaces
\begin{align} \notag
	S&:=\ker \left(-\frac{1}{2}I+\Kcal_{\Om}^{\omega_0}\right),
	\\ \notag
	S_{n,j}&:=\ker \left(-\frac{1}{2}I+P_n(T_n\Kcal_{\Om}^{\widetilde{\omega}_{n,j}})P_n\right),\quad n\ge N,\ j=1,2,\cdots,\mu.
\end{align}
We will prove that for the $L^2$-normalized computable functions in $S_{n,j}$, the components orthogonal to the exact eigenspace $S$ converge uniformly to zero as $j\to\infty$.

\begin{theorem}\label{thm:eigftn:errbdd}
	Let $\Omega$ have a $C^{1,\alpha}$ boundary for some $\alpha>0$ and set $\omega_0^2$ to be a Neumann eigenvalue of multiplicity $\mu$. If $N$ is the constant in Theorem \ref{theorem:main:second}, then
	there exists a constant $C=C(\Om,\lambda)$ such that, for all $j=1,2,\cdots,\mu$ and $n\ge N$, we have
	\beq\label{thm:eigfun:convergence} \notag
	\sup_{u\in S_{n,j},\, \|u\|_{L^2(\p\Om)}=1}\|u-\operatorname{Proj}_{S}u\|_{L^{2}(\partial\Omega)}\le C\left\|\mathcal{K}_\Omega^{\omega_0} - P_n (T_n\mathcal{K}_\Omega^{\widetilde{\omega}_{n,j}})P_n\right\|_{\mathcal{L}(L^{2}(\p\Om))},
	\eeq
	where $\operatorname{Proj}_{S}$ denotes the orthogonal projection onto the subspace $S$ of $L^{2}(\p\Om)$.
\end{theorem}

\begin{proof}
	As before, we let $A = -\frac{1}{2}I + \Kcal_{\Om}^{\omega_0}$.
	Since $-\frac{1}{2}u+P_n (T_n\mathcal{K}_\Omega^{\widetilde{\omega}_{n,j}})P_n [u] = 0$ and $A [\operatorname{Proj}_{S}u] = 0$ for $u\in S_{n,j}$ with $\|u\|_{L^{2}(\p\Om)}=1$, $n \geq N$, and $1 \leq j \leq \mu$,
	\beq\label{eq:thm:eigfun:1} \notag
	A\left[u-\operatorname{Proj}_{S}u\right]=\left(\mathcal{K}_\Omega^{\omega_0}-P_n (T_n\mathcal{K}_\Omega^{\widetilde{\omega}_{n,j}})P_n\right)[u].\eeq
	The operator $A$ has closed range because $\p\Om$ is $C^{1,\alpha}$, and therefore, 
	$$A^{-1}:\operatorname{range}A\to L^{2}(\partial\Omega)/S$$
	is bounded. We also note that 
	\beq\label{eq:def:quotnorm} \notag
	\left\|[v]\right\|_{L^{2}(\p\Om)/S}=\|v-\operatorname{Proj}_{S}v\|_{L^{2}(\p\Om)},\mbox{ for all }v\in L^{2}(\p\Om),
	\eeq
	where $[v]$ denotes the equivalence class containing $v$.
	In particular, 
	\begin{align} \notag
		\notag\|u-\operatorname{Proj}_{S}u\|_{L^{2}(\p\Om)}&=\left\|A^{-1}A\left[u-\operatorname{Proj}_{S}u\right]\right\|_{L^{2}(\p\Om)/S}\\ \notag
		\label{eq:thm:eigfun:2}&\le\left\|A^{-1}\right\|\left\|A\left[u-\operatorname{Proj}_{S}u\right]\right\|_{L^{2}(\p\Om)} \\  \notag
		& \le\left\|A^{-1}\right\|\left\|\mathcal{K}_\Omega^{\omega_0}-P_n (T_n\mathcal{K}_\Omega^{\widetilde{\omega}_{n,j}})P_n\right\|_{\mathcal{L}(L^{2}(\p\Om))}.
	\end{align}
\end{proof}

\begin{cor}\label{cor:eigftn:conv}
	For $\Omega$ with a $C^{1,\alpha}$ boundary for some $\alpha\in(0,1)$, there exists a constant $C=C(\Om,\lambda)$ such that, for all $j=1,2,\cdots,\mu$, 
	\beq\label{thm:eigfun:convergence:cor} \notag
	\lim_{n\to\infty}\, \sup_{u\in S_{n,j},\,\|u\|_{L^2(\p\Om)}=1}\|u-\operatorname{Proj}_{S}u\|_{L^{2}(\partial\Omega)}=0.
	\eeq
	
\end{cor}

\begin{proof}
	Since
	$$\left\|\mathcal{K}_\Omega^{\omega_0} - P_n (T_n\mathcal{K}_\Omega^{\widetilde{\omega}_{n,j}})P_n\right\|\le \left\|\mathcal{K}_\Omega^{\omega_0} - P_n (T_n\mathcal{K}_\Omega^{\omega_0})P_n\right\|+\left\|T_n\mathcal{K}_\Omega^{\omega_0}-T_n\mathcal{K}_\Omega^{\widetilde{\omega}_{n,j}}\right\|,$$
	Lemma \ref{Kconvergence} and part (b) of Theorem \ref{theorem:main:second} yield the desired result.
\end{proof}

\begin{remark}
	In fact, using Theorem \ref{theorem:main:smooth} and part (c) of Theorem \ref{theorem:main:second}, one can prove the convergence rate of $n^{3/2}$. 
\end{remark}


\section{Shape derivative for Eigenvalues}
\subsection{Shape Deformation using Laurent Series of Conformal Mapping}

Recall that every bounded simply connected domain $\Om$ admits a conformal mapping $\Psi_\Omega$ from $\mathbb{C}\backslash\overline{\gamma\mathbb{D}}$ onto $\mathbb{C}\backslash\overline{\Om}$ as in \eqnref{eq:Psi:laurent}. 
%
We have the following local admissibility of deformation of each Laurent series coefficient of $\Psi_\Om$:

\begin{theorem}
	Assume that $\Psi_\Omega$ admits a conformal extension to $\mathbb{C}\backslash\overline{\gamma_2\mathbb{D}}$ for some $0< \gamma_2 < \gamma$.
	Fix any $j\in\mathbb{N}$. Then there is $\varepsilon_0>0$ such that for all $\varepsilon\in\mathbb{C}$ satisfying $0\leq |\varepsilon|<\varepsilon_0$, there exists a bounded simply connected domain $\Om[k,\varepsilon]$ with an analytic boundary that has the conformal radius $\gamma$ and the exterior conformal mapping 
	\beq\label{def:Omj}
	\Psi_{\Om[k,\varepsilon]}(z)=\Psi_\Omega(z) + \varepsilon z^{-k},\quad|z|>\gamma.
	\eeq
\end{theorem}

\begin{proof}
	Fix $s_1$ and $s_2$ such that  $\gamma_2<s_2<s_1<\gamma$ and set
	$$\varepsilon_0=\frac{1}{2}{\gamma_2}^{k}\min\{|\Psi_\Om(z_1)-\Psi_\Om(z_2)|\,:\,|z_1|\ge s_1\mbox{ and }|z_2|= s_2\}>0.$$
	Let $\varepsilon\in\mathbb{C}$ satisfy $|\varepsilon|<\varepsilon_0$. We will prove that $\Psi_\Om(z)+\varepsilon z^{-k}$ is injective for $|z|>s_1$.
	
	Fix $z_0$ with $|z_0|>s_1$ and define for $\varepsilon\in\mathbb{C}$, $$f(z):=\Psi_\Om(z^{-1})-\Psi_\Om(z_0),\quad g_\varepsilon(z):=\varepsilon z^k-\varepsilon z_0^{-k},\quad 0<|z|<\gamma_2^{-1}.$$
	We remark that the only zero of $f(z)$ in $\gamma_2^{-1}\mathbb{D}$ is $z_0^{-1}$, and therefore it suffices to prove that the only zero of $f(z)+g_\varepsilon(z)$ in $\gamma_2^{-1}\mathbb{D}$ is also $z_0^{-1}$.
	From the Laurent series expansion
	$$f(z) = \frac{1}{z} + \sum_{n=0}^\infty a_n z^{n}-\Psi_\Om(z_0^{-1}),\quad 0<|z|<\gamma_2^{-1},$$
	we first note that $f(z)$ is analytic on $0<|z|<\gamma_2^{-1}$ and has a simple pole at $0$.
	Therefore, on $|z|=s_2^{-1}$, one has $$|g_\varepsilon(z)|=|\varepsilon|\left|z^{k}-z_0^{-k}\right|<\varepsilon_0\left(2\gamma_2^{-k}\right)\le |f(z)|,$$ so that by Rouche's theorem, 
	\beq\label{eq:rouche:fnfpg} \notag
	Z_f-P_f=Z_{f+g_\varepsilon}-P_{f+g_\varepsilon},
	\eeq
	where $Z_h$ and $P_h$ denote the number of zeros and poles, counting multiplicity, respectively, inside the circle $|z|=s_2^{-1}$.
	Since $f'(z_0^{-1})\ne0$, the zero $z=z_0^{-1}$ of $f$ has multiplicity one and $Z_f=1$. Using also $P_f=P_{f+g_\varepsilon}=1$, we see that $Z_{f+g_{\varepsilon}}=1$.

	We conclude that the function $\Psi_\Om(z)+\varepsilon z^{-k}$ is injective, and hence is a conformal mapping defined on $|z|>s_1$. Because the circle $|z|=\gamma$ is mapped to an analytic curve via the conformal mapping $\Psi_\Om(z)+\varepsilon z^{-k}$, the proof is complete.
	
\end{proof}

\subsection{Asymptotic analysis}\label{subsec:shapederiv}
Let $\Om[j,\varepsilon]$ be a perturbation of $\Om$ whose exterior conformal mapping is given by \eqnref{def:Omj}.
To compare $\mathcal{K}_{\Omega[k,\varepsilon]}^\omega$ and $\Kcal_\Om^\omega$ on the same function space, we define an operator
\beq\label{def:common:K2}
\begin{aligned}
	\hKcal_{\Omega[k,\varepsilon]}^\omega:\,&L^2(\p \Om)\longrightarrow L^2(\p \Om) \\
	&f\longmapsto\left(\mathcal{K}_{\Omega[k,\varepsilon]}^\omega\hat{f}\right)\circ X^{-1},\quad \hat{f}=f\circ X,
\end{aligned}
\eeq
with the change-of-variable mapping $X=\Psi_\Omega\circ{\Psi_{\Omega[k,\varepsilon]}}^{-1}:\p\Om[k,\eps]\rightarrow\p\Om$. 
One can derive the following result with some bounded integral operators $\mathcal{O}_k(\omega)$, $\widetilde{\mathcal{O}}_k(\omega)$ independent $\eps$, whose integral kernels can be explicitly stated by the asymptotic analysis similar to that in \cite{Ammari:2009:CIP}:
\begin{prop}\label{lem:ep:K}
	There exist bounded integral operators $\mathcal{O}_k(\omega)$ and $\widetilde{\mathcal{O}}_k(\omega)$ such that
	\beq\label{taylor:Khat}
	\left\|\hat{\mathcal{K}}_{\Omega[k,\varepsilon]}^\omega  -\Kcal_{\Omega}^\omega  -\varepsilon\, \mathcal{O}_k(\omega)  -\overline{\varepsilon}\,\widetilde{\mathcal{O}}_k(\omega)  \right\|_{\mathcal{L}(L^2(\p\Om))}\leq C|\varepsilon|^2 ,\eeq
	for some constant $C>0$ depending on $\Om$, $\omega>0$ and $k\in\NN$
\end{prop}

Let $\mathcal{H}=L^2(\p\Om)$. We consider an operator ${A}(\omega)$ on $\mathcal{H}$ defined by 
\beq\label{Aomega:12}
{A}(\omega)=\frac{1}{2}I+{\mathcal{K}}_\Om^\omega \quad\mbox{or}\quad -\frac{1}{2}I+{\mathcal{K}}_\Om^\omega, \text{ respectively.}
\eeq
For a positive characteristic value $\omega_0$ of ${A}(\omega)$, we then have
\beq\label{eq:mainin:pfThm1} 
{A}(\omega)^{-1}=\frac{1}{\omega-\omega_0}\mathcal{Q}+\mathcal{R}(\omega),
\eeq
where $\mathcal{Q}$ is an operator of finite rank from $\ker {A}(\omega_0)^*$ onto $\ker {A}(\omega_0)$, and $\mathcal{R}(\omega)$ is an analytic family of bounded linear operators on $\mathcal{H}$.
Note that depending on the definition of ${A}(\omega)$, $\omega_0$ is the square root of a Dirichlet or a Neumann eigenvalue of $-\Delta$.
We also set $P_0$, $P_1$, $\widetilde{P}_0$, and $\widetilde{P}_1$ to be the orthogonal projections from $\mathcal{H}$ onto 
\begin{align*}
	&P_0\mathcal{H} = {A}(\omega_0)^*\mathcal{H},\quad\qquad\widetilde{P}_0\mathcal{H} = {A}(\omega_0)\mathcal{H},\\
	&P_1\mathcal{H} = \ker{A}(\omega_0), \quad\qquad\widetilde{P}_1\mathcal{H}=\ker {A}(\omega_0)^*.
\end{align*}
By the Fredholm alternative, 
$$\mathcal{H} = P_0\mathcal{H}\oplus P_1\mathcal{H}= \widetilde{P}_0\mathcal{H}\oplus \widetilde{P}_1\mathcal{H}.$$
Note that $\dim\ker{A}(\omega_0)=\dim\ker{A}(\omega_0)^*=\mu$, where $\mu$ is the multiplicity of the Laplacian eigenvalue $\omega_0^2$. 
The operator $\mathcal{Q}$ vanishes on $\widetilde{P}_0\mathcal{H}$.
\begin{definition}\label{eta:1_mu}
	Let $\omega_0$ be the square root of a Dirichlet or a Neumann eigenvalue of $-\Delta$ on $\Om$. 
	We define an operator on $\ker{A}(\omega_0)$ as
	$$M(\omega_0,k,\theta):=\mathcal{Q}\left(e^{{\rm i}\theta}\mathcal{O}_k(\omega_0)+e^{-{\rm i}\theta}\widetilde{\mathcal{O}}_k(\omega_0)\right),\quad\mbox{for }\theta\in[0,2\pi),$$
	with $\mathcal{Q}$ given in \eqnref{Aomega:12}--\eqnref{eq:mainin:pfThm1} and $\mathcal{O}_k$, $\widetilde{\mathcal{O}}_k$ in \eqnref{taylor:Khat}. 
	Since this operator $M(\omega_0,k,\theta)$ has a $\mu\times\mu$ matrix representation, we let $\eta_1(\omega_0,k,\theta),\cdots,\eta_\mu(\omega_0,k,\theta)$ be the eigenvalues of $-M(\omega_0,k,\theta)$. The largest Jordan block for $\eta_j$ is of size $\mu_j\times \mu_j$ for $j=1,\dots,\mu$.
\end{definition}


\subsubsection{Relevant Lemmas to prove Theorem \ref{thm:shapeder:Dirichlet:multipleeig}}

\begin{lemma}[\cite{Yajima:2005:DES}]\label{lemma:blockoper:invert}
	Consider a bounded linear operator
	$L = \begin{pmatrix}
		L_{00} & L_{01} \\
		L_{10} & L_{11}
	\end{pmatrix}$
	defined on a vector space $X=X_0\oplus X_1$, where $L_{00} \in \mathcal{L}(X_0)$ is invertible. Then $L$ is invertible if and only if $E=L_{11}-L_{10}L_{00}^{-1}L_{01}$ is invertible and 
	$$L^{-1} = \begin{pmatrix}
		L_{00}^{-1} + L_{00}^{-1}L_{01}E^{-1}L_{10}L_{00}^{-1} & -L_{00}^{-1}L_{01}E^{-1} \\
		-E^{-1} L_{10}L_{00}^{-1} & E^{-1}
	\end{pmatrix}.$$
\end{lemma}
\begin{lemma}\label{lemma:inverse:norm}
	Let $\lambda$ be the eigenvalue corresponding to the largest Jordan submatrix of size $m \times m$ of an $n\times n$ matrix $M$. Then there exist $\varepsilon_0>0$ and positive constants $c$ and $C$ such that for all $\varepsilon\in\mathbb{C}$ satisfying $0<|\varepsilon|<\varepsilon_0$,
	\beq\label{eq:opnorm:inversemat} \notag
	c |\varepsilon|^{-m} \le\left\|\left(M-(\lambda+\varepsilon)I\right)^{-1}\right\|\le C |\varepsilon|^{-m}.
	\eeq
\end{lemma}
\begin{proof}
	Let $\lambda I + N$ be a Jordan block of size $m\times m$ of $M$, where $[N]_{ij} = \delta_{i+1,j}$.
	Then $-\varepsilon I + N$ is a Jordan block of $M-(\lambda+\varepsilon)I$ and since $N^m=0$,
	$$(-\varepsilon I + N)^{-1} = (-\varepsilon)^{-1} \sum_{j=0}^{m-1} (\varepsilon^{-1} N)^{j}.$$
	Due to the fact that $\|N\|=1$, we then have
	\beq\label{eq:uppbd:jordanblock} \notag
	\|(-\varepsilon I + N)^{-1}\|\le\sum_{j=1}^m |\varepsilon|^{-j} \le C|\varepsilon|^{-m},
	\eeq
	for some constant $C$.
	
	Next, we set $e_k$ to be the $k$-th column of $I_m$. Since $N^je_m = e_{m-j}$ for $1\le j\le m-1$, 
	\beq\label{eq:lwrbd:jordanblock} \notag
	\|(-\varepsilon I + N)^{-1}\| \ge \|(-\varepsilon I + N)^{-1}e_m\| = \left\|(-\varepsilon)^{-1} \sum_{j=0}^{m-1} \varepsilon^{-j} e_{m-j}\right\| = \left(\sum_{j=1}^{m}|\varepsilon|^{-2j}\right)^{1/2}\ge|\varepsilon|^{-m}.
	\eeq
	Finally, we consider the Jordan decomposition $M=S\operatorname{diag}(B_1,B_2,\cdots,B_n)S^{-1}$, where $S$ is an invertible matrix and the $B_j's$ are Jordan blocks with diagonal entries $\lambda_j$.
	Let $0<\varepsilon_0<1$ satisfy $\sigma(M)\cap B(\lambda,\sqrt{\varepsilon_0})=\{\lambda\}$.
	For $\varepsilon\in\mathbb{C}$ such that $0<|\varepsilon|<\varepsilon_0$, we set $\widetilde{B}_j:=B_j-(\lambda+\varepsilon)I$.
	Then 
	$$\|(M-(\lambda+\varepsilon)I)^{-1}\| = \left\|S\left(\operatorname{diag}(\widetilde{B}_1^{-1},\widetilde{B}_2^{-1},\cdots,\widetilde{B}_n^{-1})\right)S^{-1}\right\|\le\|S\|\left(\sum_{j=1}^n\|\widetilde{B}_j^{-1}\|\right)\|S^{-1}\|,$$
	and this completes the proof.
\end{proof}

\subsubsection{Proof of Theorem \ref{thm:shapeder:Dirichlet:multipleeig}}

Note that
$$\left(\Kcal_\Om^\omega-\Kcal_\Om^{\omega_0}\right)[f](x)=\int_{\p\Om}\frac{\p (\Gamma^\omega-\Gamma^{\omega_0})(x-y)}{\p \nu_y}f(y)\,d\sigma(y),$$
where $\Gamma^\omega-\Gamma^{\omega_0}$ is regular with respect to $\omega$ as well as $x$. 
One can show that
\beq\label{HKexp:2}
\hKcal_{\Omega[k,\eps]}^\omega -\hKcal_{\Omega[k,\eps]}^{\omega_0}=\Kcal_\Om^\omega-\Kcal_\Om^{\omega_0}+T_0(\eps,\omega),\quad\mbox{where }\|T_0(\eps,\omega)\|\leq O(|\eps||\omega-\omega_0|).
\eeq
By \eqnref{taylor:Khat} and \eqnref{HKexp:2},
\beq\label{hatA:esti}
\hat{\mathcal{K}}_{\Omega[k,\varepsilon]}^\omega - {\mathcal{K}}_{\Omega}^\omega = \varepsilon\mathcal{O}_{k}(\omega_0)+\overline{\varepsilon}\,\widetilde{\mathcal{O}}_{k}({\omega_0})+ T(\varepsilon,\omega),\quad\mbox{where } \|T(\varepsilon,\omega)\|=O(|\varepsilon|^2+|\eps||\omega-\omega_0)|).
\eeq

\begin{proof}[Proof of Theorem \ref{thm:shapeder:Dirichlet:multipleeig}]
	We set $\Delta\hat{A}[k,\varepsilon]:=\hat{\mathcal{K}}_{\Omega[k,\varepsilon]}^\omega - {\mathcal{K}}_{\Omega}^\omega$.
	Let $U$ be the unitary operator satisfying
	$$U\widetilde{P}_1\mathcal{H} = P_1\mathcal{H}\quad\mbox{and}\quad U\widetilde{P}_0\mathcal{H}=P_0\mathcal{H}.$$
	We apply Lemma \ref{lemma:blockoper:invert} to 
	$$X_i = P_i\mathcal{H},\quad L = U \left({A}(\omega)+\Delta\hat{A}[k,\varepsilon]\right),\quad L_{im}=P_i L P_m\quad\mbox{for }i,m=0,1.$$
	Note that $L=U\left(\frac{1}{2}I+\Kcal_\Om^\omega\right)$ for the Dirichlet eigenvalue case and $L=U\left(-\frac{1}{2}I+\Kcal_\Om^\omega\right)$ for the Neumann eigenvalue case.
	Applying the expansion ${A}(\omega)=A_0 + (\omega-\omega_0)A_1(\omega)$ with $A_0=A(\omega_0)$, we have
	\begin{align*}
		&L_{00}=UA_0 + (\omega-\omega_0)P_0 U A_1(\omega) P_0 + P_0 U \Delta\hat{A}[k,\varepsilon] P_0,\\
		&L_{im}=(\omega-\omega_0)P_i UA_1(\omega) P_m + P_i U \Delta\hat{A}[k,\varepsilon] P_m,\quad\mbox{for }(i,m) = (0,1),\,(1,0),\,(1,1).
	\end{align*}
	The operator $UA_0$ is invertible on $P_0\mathcal{H}$. 
	Furthermore, by \eqnref{hatA:esti},
	$$\Delta\hat{A}[k,\varepsilon] = \varepsilon\mathcal{O}_k(\omega_0)+\overline{\varepsilon}\,\widetilde{\mathcal{O}}_k(\omega_0)+ T(\varepsilon,\omega),\quad\mbox{where}\quad \|T(\varepsilon)\|=O(|\varepsilon|^2+|\eps||\omega-\omega_0)|).$$
	Hence, there exists some $\varepsilon_0>0$ such that for all $\varepsilon,\omega\in\mathbb{C}$ satisfying $|\omega-\omega_0|<\varepsilon_0$ and $|\varepsilon|<\varepsilon_0$, the operator $L_{00}$ is invertible.
	
	Next, we address the invertibility of $E=L_{11}-L_{10}L_{00}^{-1}L_{01}$. For simplicity, we write $M(\omega_0,k,\arg(\eps))=M(\arg(\eps))$ and $\eta_j(\omega_0,k,\arg(\omega))=\eta_j(\arg(\omega))$. 
	It can be shown using $(P_1UA_1(\omega_0)P_1)^{-1}=\mathcal{Q}U^{-1}$ that 
	\begin{align*}
		&L_{11}=L_{11}^{(1)}+L_{11}^{(2)}, \text{ where} \\ 
		&L_{11}^{(1)}:=|\varepsilon| \left(P_1UA_1(\omega_0)P_1\right)\left(\frac{\omega-\omega_0}{|\varepsilon|}I +  M\left(\arg(\varepsilon)\right)\right),\\
		&L_{11}^{(2)}:=  P_1 U \left((\omega-\omega_0)(A_1(\omega)-A_1(\omega_0))+T(\varepsilon,\omega)\right) P_1.
	\end{align*}
	From Lemma \ref{lemma:inverse:norm}, there exist positive constants $\varepsilon_1$ and $C$ such that for all $\omega\in\mathbb{C}$ satisfying $0<|\omega-\omega_0-|\varepsilon|\eta_j(\arg(\varepsilon))|<\varepsilon_1|\varepsilon|$,
	$$\left\|\left(\frac{\omega-\omega_0}{|\varepsilon|}I + M\left(\arg(\varepsilon)\right)\right)^{-1}\right\|\le C \left|\frac{\omega-\omega_0}{|\varepsilon|}-\eta_j(\arg(\varepsilon))\right|^{-\mu_n}.$$
	We choose $C$ so that $\|T(\varepsilon,\omega)\|\le C(|\varepsilon|^2+|\eps||\omega-\omega_0|)$, $\|L_{10}\|\le C|\varepsilon|$, $\|L_{01}\|\le C|\varepsilon|$, and $\|L_{00}^{-1}\|\le C$ also hold.
	Because $A_1(\omega)$ is analytic in $\omega$, there also exists a constant $C'$ such that 
	$$\|A_1(\omega)-A_1(\omega_0)\|\le C'|\omega-\omega_0|, \quad\omega\in B(\omega_0,|\omega_0|/2).$$
	For $\omega\in\mathbb{C}$ on the circle
	\beq\label{eq:circle:firstorder}
	\left|\omega-\omega_0-|\varepsilon|\eta_j(\arg(\varepsilon))\right|=C_0|\varepsilon|^{1+\frac{1}{\mu_j}},
	\eeq
	where $C_0$ is a constant to be determined later, we have
	\begin{align*}
		\|L_{11}^{-1}\| & = \left\|\big(L_{11}^{(1)}\big)^{-1}\sum_{m=0}^\infty\left(-\big(L_{11}^{(1)}\big)^{-1}L_{11}^{(2)}\right)^m\right\|\\
		&\le CC_0^{-\mu_j}|\varepsilon|^{-2}\sum_{m=0}^\infty \left(CC_0^{-\mu_j}|\varepsilon|^{-2} \left(C'\left(|\varepsilon|\eta_j(\arg(\varepsilon))+C_0|\varepsilon|^{1+\frac{1}{\mu_j}}\right)^2+C|\varepsilon|^2\right)\right)^m.
	\end{align*}
	For $C_0$ large enough, $\|L_{11}^{-1}\|\le c|\varepsilon|^{-2}$ so that
	$$\|E^{-1}\|=\left\|L_{11}^{-1}\sum_{m=0}^\infty(L_{11}^{-1}L_{10}L_{00}^{-1}L_{01})^m \right\|\le c|\varepsilon|^{-2},$$
	where $c\to0$ as $C_0\to\infty$.
	We choose a sufficently small $c$ so that 
	$$\|E^{-1} \|\left\|L_{11}^{(2)}-L_{10}L_{00}^{-1}L_{01}\right\|<1,$$
	on the circle given by \eqnref{eq:circle:firstorder}.
	
	We now apply the generalized Rouche's theorem to the operators $E$ and $L_{11}^{(1)}$, where the characteristic values of $L_{11}^{(1)}$ are $\omega=\omega_0+|\varepsilon|\eta_j(\arg(\varepsilon))$. We conclude that there exists a characteristic value of $E$, and thus of $L$, inside the circle \eqnref{eq:circle:firstorder}.

	\smallskip
	The constant $C$ in \eqnref{eq:asympt:Dirichlet:simplified} is independent of both $|\varepsilon|$ and $\arg(\varepsilon)$.
	The independence with respect to $\arg(\varepsilon)$ follows from the fact that each $\eta_n(\theta)$ is bounded uniformly in $\theta$, which can be seen as follows:
	since
	$$|[M(\omega_0,k,\theta)]_{mn}|\le |[\mathcal{Q}\mathcal{O}_j]_{mn}| + |[\mathcal{Q}\widetilde{\mathcal{O}}_j]_{mn}|,\quad 1\le m,n\le \mu,$$
	the coefficients of the characteristic polynomial
	$$\det(\eta I - [M(\omega_0,k,\theta)])=\eta^\mu + \sum_{m=0}^{\mu-1}c_m \eta^m$$
	satisfy $|c_m|\le C$ for a constant $C$ independent of $\theta$.
	Therefore, the roots $\eta_n(\theta)$ of the characteristic polynomial are bounded by 
	$$|\eta_n(\theta)|\le\max\left\{1,\sum_{m=0}^{\mu-1}|c_m|\right\}\le\max\{1,\mu C\}.$$
\end{proof}

\bigskip


\end{document}